\documentclass[11pt, reqno]{amsart}

\usepackage[OT2, T1]{fontenc}

\usepackage[usenames,dvipsnames]{color}

\usepackage{url}
\usepackage{amsmath}
\usepackage{array}
\usepackage{amsfonts}
\usepackage{amsfonts}
\usepackage{amssymb}
\usepackage{amsthm}
\definecolor{link}{RGB}{11,0,128}
\usepackage[colorlinks=true,citecolor=NavyBlue,linkcolor=Bittersweet,urlcolor=link]{hyperref}
\usepackage{colonequals}	% for nice :=
\usepackage{color}
\usepackage{mathabx}		% for \widec
\usepackage{mathrsfs}
\usepackage{amscd}
\usepackage[all,cmtip]{xy}	% for commutative diagrams
\usepackage{sseq}			% for spectral sequences
\usepackage{verbatim}
\usepackage{parskip}
\usepackage{microtype}
\usepackage[in]{fullpage}
\usepackage{graphicx}
\usepackage{mathrsfs} 			% for \mathscr (script letters)
\usepackage{bm} 				% for bold Greek letters
\usepackage{cleveref}
\usepackage{enumitem}
\usepackage{moreenum}			% for enumeration via Greek letters
\usepackage[alphabetic,lite]{amsrefs} 	% for bibliography
\usepackage{pifont}
\usepackage{stmaryrd}	% this is for \llbracket and \rrbracket
\usepackage{graphicx}
\usepackage{rotating}

\DeclareSymbolFont{cyrletters}{OT2}{wncyr}{m}{n}
\DeclareMathSymbol{\Sha}{\mathalpha}{cyrletters}{"58}

\usepackage{color}

\newcommand{\gA}{\alpha}
\newcommand{\gB}{\beta}

\newcommand{\gL}{\lambda}

\newcommand{\bC}{\mathbb{C}}

\newcommand{\bF}{\mathbb{F}}
\newcommand{\bG}{\mathbb{G}}

\newcommand{\bQ}{\mathbb{Q}}
\newcommand{\bR}{\mathbb{R}}

\newcommand{\bZ}{\mathbb{Z}}

\newcommand{\cA}{\mathcal{A}}
\newcommand{\cB}{\mathcal{B}}

\newcommand{\cG}{\mathcal{G}}
\newcommand{\cH}{\mathcal{H}}

\newcommand{\cO}{\mathcal{O}}

\newcommand{\cX}{\mathcal{X}}

\newcommand{\fo}{\mathfrak{o}}

\newcommand{\sL}{\mathscr{L}}

\newcommand{\ra}{\rightarrow}

\newcommand{\Ra}{\Rightarrow}

\newcommand{\xra}{\xrightarrow}

\newcommand{\hra}{\hookrightarrow}

\newcommand{\I}{^{\infty}}
\newcommand{\wt}{\widetilde}
\newcommand{\wh}{\widehat}
\newcommand{\eps}{\epsilon}

\newcommand{\pr}{^{\prime}}

\newcommand{\ce}{\colonequals}
\newcommand{\ov}{\overline}

\renewcommand{\b}{\textbf}

\newcommand{\Qlbar}{\ov{\bQ}_\ell}		% Algebraic closure of l-adic numbers

\newcommand{\tensor}{\otimes} 		% binary tensor product
\newcommand{\isomto}{\overset{\sim}{\longrightarrow}}

\newcommand{\nd}{{\mathrm{nd}}}		% the maximal nondivisible quotient of an abelian group
\newcommand{\fppf}{\mathrm{fppf}}		% for fppf cohomology (mainly used in subscripts)
\newcommand{\et}{\mathrm{\acute{e}t}}	% for etale cohomology (mainly used in subscripts)
\newcommand{\tors}{\mathrm{tors}}		% torsion subgroup (mainly used in subscripts)
\renewcommand{\ss}{\mathrm{ss}}		% semisimplification (mainly used in superscripts}
\newcommand{\Frob}{\mathrm{Frob}}		% Frobenius
\renewcommand{\i}{^{-1}}
\renewcommand{\th}{^{\mathrm{th}}}

\providecommand{\abs}[1]{\left\lvert#1\right\rvert}

\providecommand{\In}[1]{\left\langle#1\right\rangle}
\providecommand{\p}[1]{\left(#1\right)}

\providecommand{\f}[2]{\frac{#1}{#2}}
\DeclareMathOperator{\Ker}{Ker}			% Kernel
\DeclareMathOperator{\Coker}{Coker}		% Cokernel
\DeclareMathOperator{\im}{Im}			% Imaginary part
\DeclareMathOperator{\Spec}{Spec}		% Spectrum of a ring
\DeclareMathOperator{\Hom}{Hom}			% Set of arrows between two object

\DeclareMathOperator{\Char}{char}		% Characteristic of a field
\DeclareMathOperator{\id}{id}			% identity
\DeclareMathOperator{\Ext}{Ext}			% Derived functors of Hom
\DeclareMathOperator{\loc}{loc}		% 
\DeclareMathOperator{\Br}{Br}		% Brauer group
\DeclareMathOperator{\Gal}{Gal}	% Galois group
\DeclareMathOperator{\tr}{tr}		% Truncation of a simplicial object or a trace
\DeclareMathOperator{\Tr}{Tr}		% Trace
\DeclareMathOperator{\inv}{inv}	% The invariant map from local class field theory
\DeclareMathOperator{\ord}{ord}	% order
\DeclareMathOperator{\Ind}{Ind}		% Induced representation
\DeclareMathOperator{\Res}{Res}		% Restriction of the representation
\DeclareMathOperator{\rk}{rk}		% rank
\DeclareMathOperator{\Sel}{Sel}		% Selmer group
\newcommand{\x}{\text}
\newcommand{\qq}{\quad\quad}
\newcommand{\qqq}{\quad\quad\quad}

\newcommand{\tst}{\textstyle}

\DeclareMathOperator{\Sw}{Sw}		% Swan conductor
\DeclareMathOperator{\NS}{NS}		% Neron--Severi group
\newcommand{\ba}{\begin{aligned}}
\newcommand{\ea}{\end{aligned}}
\newcommand{\be}{\begin{equation}}
\newcommand{\ee}{\end{equation}}
\newcommand{\pf}{\begin{proof}}
\newcommand{\bpf}{\begin{proof}}
\newcommand{\epf}{\end{proof}}
\newcommand{\bthm}{\begin{thm}}
\newcommand{\ethm}{\end{thm}}
\newcommand{\bprop}{\begin{prop}}
\newcommand{\eprop}{\end{prop}}
\newcommand{\bcor}{\begin{cor}}
\newcommand{\ecor}{\end{cor}}
\newcommand{\brem}{\begin{rem}}
\newcommand{\erem}{\end{rem}}
\newcommand{\brems}{\begin{rems} \hfill \begin{enumerate}[label=\b{\thesubsection.},ref=\thesubsection]}

\newcommand{\erems}{\end{enumerate} \end{rems}}
\newcommand{\blem}{\begin{lemma}}
\newcommand{\elem}{\end{lemma}}
\newcommand{\bconj}{\begin{conj}}
\newcommand{\econj}{\end{conj}}
\newcommand{\bprob}{\begin{Problem}}
\newcommand{\eprob}{\end{Problem}}
\newcommand{\bq}{\begin{q}}
\newcommand{\eq}{\end{q}}
\newcommand{\benum}{\begin{enumerate}[label={{\upshape(\alph*)}}]}
\newcommand{\eenum}{\end{enumerate}}
\newcommand{\benumr}{\begin{enumerate}[label={{\upshape(\roman*)}}]}
\newcommand{\bc}{\begin{comment}}
\newcommand{\ec}{\end{comment}}
\newcommand{\beg}{\begin{eg}}
\newcommand{\eeg}{\end{eg}}
\newcommand{\lab}{\label}

\theoremstyle{plain}
\newtheorem{thm}[subsection]{Theorem}
\Crefname{thm}{Theorem}{Theorems}

\Crefname{rethm}{Theorem}{Theorem}
\newtheorem{prop}[subsection]{Proposition}
\Crefname{prop}{Proposition}{Propositions}
\newtheorem{q}[subsection]{Question}
\Crefname{q}{Question}{Questions}
\Crefname{eg}{Example}{Examples}
\newtheorem{Problem}[subsection]{Problem}
\Crefname{Problem}{Problem}{Problems}
\newtheorem{conj}[subsection]{Conjecture}
\Crefname{conj}{Conjecture}{Conjectures}
\newtheorem{cor}[subsection]{Corollary}
\Crefname{cor}{Corollary}{Corollaries}
\newtheorem{lem}[equation]{Lemma}
\Crefname{lem}{Lemma}{Lemmas}

\newtheorem{lemma}[subsection]{Lemma}

\theoremstyle{remark}

\Crefname{claim}{Claim}{Claims}

\theoremstyle{definition}
\newtheorem{eg}[subsection]{Example}
\newtheorem{rem}[subsection]{Remark}
\Crefname{rem}{Remark}{Remarks}
\newtheorem*{rems}{Remarks}

\newtheoremstyle{subsection-tweak}
   {11pt}
   {3pt}%
   {}
   {}%
   {\bfseries}
   {}%
   {.5em}
   {\thmnumber{\@{#1}{}\@{#2}.}%
    \thmnote{~{\bfseries#3.}}}

\Crefname{innercustomconj}{Conjecture}{Conjecture}

\theoremstyle{subsection-tweak}
\newtheorem{pp}[subsection]{}
\newcommand{\bpp}{\begin{pp}}
\newcommand{\epp}{\end{pp}}

\numberwithin{equation}{subsection}

\usepackage[foot]{amsaddr}

\begin{document}
\author{K\k{e}stutis \v{C}esnavi\v{c}ius}
\title{\resizebox{\textwidth}{!}{\mbox{The $\ell$-parity conjecture over the constant quadratic extension}}}
\date{\today}
\address{Department of Mathematics, University of California, Berkeley, CA 94720-3840, USA}
\email{kestutis@berkeley.edu}

\begin{abstract} 
For a prime $\ell$ and an abelian variety $A$ over a global field $K$, the $\ell$-parity conjecture predicts that, in accordance with the ideas of Birch and Swinnerton-Dyer, the $\bZ_\ell$-corank of the $\ell\I$-Selmer group and the analytic rank agree modulo $2$. Assuming that $\Char K > 0$, we prove that the $\ell$-parity conjecture holds for the base change of $A$ to the constant quadratic extension if $\ell$ is odd, coprime to $\Char K$, and does not divide the degree of every polarization of $A$. The techniques involved in the proof include the \'{e}tale cohomological interpretation of Selmer groups, the Grothendieck--Ogg--Shafarevich formula, and the study of the behavior of local root numbers in unramified extensions.
\end{abstract}

\maketitle

% 2013.12.9 -- 2014.2.9

\section{Introduction}

\bpp[The $\ell$-parity conjecture] \lab{l-par}
The Birch and Swinnerton-Dyer conjecture (BSD) predicts that the completed $L$-function $L(A, s)$ of an abelian variety $A$ over a global field $K$ extends meromorphically to the whole complex plane, vanishes to the order $\rk A$ at $s = 1$, and satisfies the functional equation
\be\lab{fe}
L(A, 2 - s) = w(A) L(A, s),
\ee
where $\rk A$ and $w(A)$ are the Mordell--Weil rank and the global root number of $A$. The vanishing assertion combines with \eqref{fe} to give ``BSD modulo $2$,'' namely, the parity conjecture:
\[
(-1)^{\rk A} \overset{?}{=} w(A).
\]
Selmer groups tend to be easier to study than Mordell--Weil groups, so, fixing a prime $\ell$, one sets
\[
\rk_\ell A \ce \dim_{\bQ_\ell} \Hom(\Sel_{\ell\I} A, \bQ_\ell/\bZ_\ell)\tensor_{\bZ_\ell} \bQ_\ell,
\]
where  
\[
\Sel_{\ell\I} A \ce \varinjlim \Sel_{\ell^n} A
\]
is the $\ell\I$-Selmer group, notes that the conjectured finiteness of the Shafarevich--Tate group $\Sha(A)$ implies $\rk_\ell A = \rk A$, and instead of the parity conjecture considers the $\ell$-parity conjecture:
\be\lab{l-par-conj}
(-1)^{\rk_\ell A} \overset{?}{=} w(A).
\ee
\epp

\bpp[The status of the number field case] \lab{prev}
Over a number field $K$, in the elliptic curve case, the $\ell$-parity conjecture was proved for $K = \bQ$ in \cite{DD10}*{Thm.~1.4} (see \cite{Mon96}, \cite{Nek06}*{\S0.17}, \cite{Kim07} for some preceding work), for totally real $K$ excluding some cases of potential complex multiplication with $\ell = 2$ in \cite{Nek13}*{Thm.~A}, \cite{Nek15}*{5.12}, and \cite{Nek16}*{Thm.~E}, and for curves with an $\ell$-isogeny over a general $K$ in \cite{Ces16d}*{Thm.~1.4} (building on \cite{DD08a}*{Thm.~2}, \cite{DD11}*{Cor.~5.8}, and \cite{CFKS10}*{proof of Thm.~2.1}). The higher dimensional case in the presence of a suitable isogeny was addressed in \cite{CFKS10}*{Thm.~2.1}.
\epp

\bpp[The status of the function field case] \lab{ff-case}
Over a global field $K$ of positive characteristic $p$, the elliptic curve case of the $p$-parity conjecture was proved in \cite{TW11}*{Thm.~1} under the assumption that $p > 3$. The function field case of the $\ell$-parity conjecture was subsequently settled in full in \cite{TY14}*{Thm.~1.1} (including the case $\ell = p$).
 
The main goal of this paper is to present another approach to the $\ell$-parity conjecture in the positive characteristic case. Even though our techniques only reprove the special case stated in \Cref{main}, they lead to intermediate results that have already been useful in other contexts---for instance, in the proof of the Kramer--Tunnell formula for certain classes of hyperelliptic Jacobians presented in the PhD thesis \cite{Mor15} of Adam Morgan.

%Assume that $\Char K = p$ with $p > 0$. Then the meromorphic continuation of $L(A, s)$ is known \cite{Sch82}*{Bemerkung 2) on p.~497}, as is the functional equation \eqref{fe} \cite{Del73}*{8.10 and 9.3}.\footnote{Note the misprints: $pt\i$ should be $(pt)\i$ in \cite{Del73}*{p.~4, 7.11 (iii), 7.11.8 and the following equation, 9.3.1, 9.9.1, 9.9.2}. Also note the need of including in $L(A, s)$ the factors contributed by the conductor to get the claimed form of \eqref{fe}.} %\kestutis{Also note the need to include the factors contributed by the discriminant, too.}
%Moreover, as observed in \cite{TW11}*{Prop.~3}, \cite{KT03} proves $\rk_\ell A \le \ord_{s = 1} L(A, s)$ for every $\ell$ (for $\ell \neq p$, already \cite{Sch82}*{Lemma~2~i)} suffices). Also, due to \cite{KT03} (which extends the results of \cite{Sch82}), the equality in $\rk A \le \rk_\ell A$ for a single $\ell$ implies the equality for every $\ell$ along with the full Birch and Swinnerton-Dyer conjecture for $A$. Consequently, $\rk_\ell A = \ord_{s = 1} L(A, s)$ if $\ord_{s = 1} L(A, s) \le 1$, so the $\ell$-parity conjecture holds in this case. Knowing it in general amounts to knowing that the nonnegative integer $\ord_{s = 1} L(A, s) - \rk_\ell A$, which is expected to be $0$, is at least even. We prove this in most cases over the constant quadratic extension:
\epp

\bthm[\Cref{main-pf}] \lab{main}
Let $K$ be a global field of positive characteristic, let $\bF_q$ be its field of constants, let $\ell \nmid q$ be a prime, and let $A$ be an abelian variety over $K$. Suppose that $A$ has a polarization of degree prime to $\ell$ {\upshape(}e.g.,~a principal polarization{\upshape)} and, if $\ell = 2$, that the orders of the component groups of the reductions of $A_{K\ov{\bF}_q}$ are odd. Then the $\ell$-parity conjecture holds for~$A_{K\bF_{q^2}}$.
\ethm

\brem
\Cref{l=2} isolates the difficulty in removing the additional component group condition in the case $\ell = 2$.
\erem

One of the key inputs to the proof of \Cref{main} is the following purely local result that allows us to control the root numbers over an unramified quadratic extension.

\bthm[\Cref{st-root}]\lab{st-root-ann}
For an abelian variety $B$ over a nonarchimedean local field $k$, let $B_{k_n}$ and $a(B)$ be its base change to a degree $n$ unramified extension and conductor exponent, respectively. The local root number satisfies
\[
w(B_{k_n}) = \begin{cases} w(B),\quad \quad\  \text{ if $n$ is odd,} \\  (-1)^{a(B)},\quad \text{ if $n$ is even.}  \end{cases}
\]
\ethm

\Cref{st-root-ann} leads to the following extension of a theorem of Kisilevsky, \cite{Kis04}*{Thm.~1}, who treated the case when $K = \bQ$ and $B$ and $B'$ are elliptic curves.

\bthm
Let $B$ and $B'$ be abelian varieties over a global field $K$. If $w(B_{L}) = w(B'_L)$ for every separable quadratic extension $L/K$, then the conductor ideals of $B$ and $B'$ are equal up to square factors.
\ethm

\bpf
Let $\Sigma$ be the set of those finite places of $K$ at which $B$ or $B'$ has bad reduction. For a $v \in \Sigma$, let $L/K$ be a separable quadratic extension that is inert at $v$, split at every $v' \in \Sigma$ with $v' \neq v$, and split at every infinite place. For such an $L$ we have
\[
w(B_L) = w(B_{L_v}) \qq \text{and} \qq w(B'_L) = w(B'_{L_v}),
\]
so \Cref{st-root-ann} and the $w(B_L) = w(B'_L)$ assumption show that up to squares the conductor ideals of $B$ and $B'$ have the same factor at $v$.
\epf

\bpp[An overview of the proof of \Cref{main}] \lab{outline}
The main idea of the proof is to interpret the $\ell$-Selmer group as an \'{e}tale cohomology group following \cite{Ces16c}, to use the Grothendieck--Ogg--Shafarevich formula together with the Hochschild--Serre spectral sequence to express the $\ell$-Selmer parity as a sum of local terms, and then to compare place by place with the expression of the global root number as the product of local root numbers. The argument is carried out in \S\ref{final} and rests on the following additional inputs:
\begin{itemize}
\item
Due to possibly nontrivial Galois action, the Hochschild--Serre spectral sequence relates the \'{e}tale cohomological $\ell$-Selmer group to the Grothendieck--Ogg--Shafarevich formula only after a large constant extension. In order to descend the $\ell$-parity conclusion to $K\bF_{q^2}$, we use the self-duality of Galois representations furnished by $\ell\I$-Selmer groups. This self-duality is the subject of \S\ref{self-dual}.

\item
For the main idea to be relevant, in \S\ref{sel-sel} we implicitly recast the $\ell$-parity conjecture in terms of $\ell$-Selmer rather than $\ell\I$-Selmer groups. Special care is needed if $\ell = 2$, since Shafarevich--Tate groups need not be of square order even when they are finite and $A$ is principally polarized. This is one of the places of the overall argument where the polarization assumption of \Cref{main} comes in (the polarization is also used in \S\S\ref{ade}--\ref{final} to get an isomorphism $A[l]\, \simeq A^\vee[l]$ and to apply suitable results from \cite{PR12}).

\item
The comparison of the Grothendieck--Ogg--Shafarevich local terms and root numbers is possible due to \Cref{st-root-ann}, which, along with related local results, is treated in \S\ref{app}.

\item
In the presence of local Tamagawa factors that are divisible by $\ell$, the $\ell$-Selmer group may differ from its \'{e}tale cohomological counterpart. Arithmetic duality results proved in \S\S\ref{ad-pre}--\ref{ade} control this difference modulo $2$ through \Cref{dual}. The $\ell = 2$ case again leads to complications due to the difference between alternating and skew-symmetric pairings in characteristic $2$.
\end{itemize}
\epp

\bpp[Notation]\lab{not}
For a field $F$, an algebraic (resp.,~separable) closure is denoted by $\ov{F}$ (resp.,~by $F^s$);  when needed (e.g., for forming composita), the choice of $\ov{F}$ is made implicitly and compatibly with overfields. If $F$ is a global field and $v$ is its place, then $F_v$ denotes the corresponding completion; if $v$ is finite, then $\cO_v$ and $\bF_v$ denote the ring of integers and the residue field of $F_v$. For a prime $\ell$ and a torsion abelian group $G$, we denote by $G[\ell\I]$ (resp.,~by $G_\nd$) the subgroup consisting of all the elements of $\ell$-power order (resp.,~the quotient by the maximal divisible subgroup). In the case when $G = G[\ell\I]$, we say that $G$ is \emph{$\bZ_\ell$-cofinite} if $\Hom(G, \bQ_\ell/\bZ_\ell)$ is finitely generated as a $\bZ_\ell$-module. For a prime $\ell$ and an abelian variety $B$ over a global field, the notation $\rk B$, $w(B)$, $\rk_\ell B$, $\Sel_{\ell\I} B$, $\Sha(B)$ introduced in \S\ref{l-par} remains in place throughout the paper. If $B$ is an abelian variety over a local field, then $w(B)$ denotes the \emph{local} root number instead. The dual abelian variety is denoted by $B^\vee$. All the representations that we consider are finite dimensional.
\epp

\subsection*{Acknowledgements}
I thank the referee for a careful reading of the manuscript and for very helpful comments. I thank Tim Dokchitser, Vladimir Dokchitser, Adam Morgan, Bjorn Poonen, Padmavathi Srinivasan, Fabien Trihan, and Christian Wuthrich for helpful conversations or correspondence.

\section{Self-duality of $\ell\I$-Selmer groups} \lab{self-dual}

Let $F/K$ be a finite Galois extension of global fields, let $\ell$ be a prime number, and let $A$ be an abelian variety over $K$. The goal of this section is to prove in \Cref{DD} that if $\ell \neq \Char K$, then the finite-dimensional $\bQ_\ell$-representation
\[
\cX_\ell(A) \ce \Hom(\Sel_{\ell\I} A, \bQ_\ell/\bZ_\ell) \tensor_{\bZ_\ell} \bQ_\ell
\]
of $\Gal(F/K)$ is isomorphic to its dual and to deduce in \Cref{odd-gal} that if, in addition, the degree of $F/K$ is odd, then the $\ell$-parity conjecture holds for $A$ if and only if it holds for $A_F$. The self-duality is expected: if $\Sha(A_F)[\ell\I]$ is finite as is conjectured, then $\cX_\ell(A)$ is $\Gal(F/K)$-isomorphic to $A(F) \tensor_{\bZ} \bQ_\ell$, which is self-dual. The utility of \Cref{DD} lies in bypassing assumptions on $\Sha$.

In the number field case, proofs of these results were given by T.~and V.~Dokchitser, see \cite{DD09b}*{Prop.~A.2} for a summary. In general, their arguments require only mild modifications and are explained in the proof of \Cref{DD} below. In \cite{TW11}*{proof of Prop.~4}, Trihan and Wuthrich have also observed that extensions of the sort presented here are possible, but it seems worthwhile to indicate the necessary changes to the proofs.

A key input to the proof of the  promised self-duality of $\ell^\infty$-Selmer groups is the following lemma, whose proof is based on the Selmer group analogue of the isogeny invariance of BSD quotients.

\blem \lab{Mil-lem}
Let $\phi \colon X \ra Y$ be an isogeny between abelian varieties over a global field $K$ and set
\[
Q(\phi) \ce \prod_{\ell \mid \deg \phi} \# \Coker\p{ \f{\Hom(\Sel_{\ell\I} Y, \bQ_\ell/\bZ_\ell)}{\Hom(\Sel_{\ell\I} Y, \bQ_\ell/\bZ_\ell)_\tors} \xra{\phi^*} \f{\Hom(\Sel_{\ell\I} X, \bQ_\ell/\bZ_\ell)}{\Hom(\Sel_{\ell\I} X, \bQ_\ell/\bZ_\ell)_\tors} }, 
\]
and likewise for $Q(\phi^\vee)$. If $X = Y$ and $\deg \phi$ is prime to $\Char K$, then $Q(\phi) = Q(\phi^\vee)$.
\elem

\bpf
The proof of the number field case \cite{DD10}*{Thm.~4.3} continues to work (as in \emph{loc.~cit.},~one bases the argument on \cite{Mil06}*{I.(7.3.1)}). 
\epf

\bthm \lab{DD}
For a finite Galois extension $F/K$ of global fields, a prime $\ell \neq \Char K$, and an abelian variety $A$ over $K$, the $\bQ_\ell$-representation $\cX_\ell(A_{F})$ of $G \ce \Gal(F/K)$ is self-dual.
\ethm

\bpf
We model the argument on the proof, given by T.~and V.~Dokchitser in \cite{DD09c}*{Thm.~2.1}, of the corresponding statement in the number field case. To begin with, since $A$ and $A^\vee$ are isogenous, $\cX_\ell(A_F) \simeq \cX_\ell(A^\vee_F)$ as $\bQ_\ell[G]$-modules, so the Zarhin trick allows us to assume that $A$ is principally polarized. 

We will apply \Cref{Mil-lem} with $X = Y = \Res_{F/K}(A_F)$, which, due to the functoriality of the restriction of scalars $\Res_{F/K}$, comes equipped with a $G$-action. By the proof of \cite{DD09c}*{Lemma 2.4} (originally given in the number field case), there is an isomorphism 
\[
\cX_\ell(\Res_{F/K}(A_F)) \cong \cX_\ell(A_F) \qq \text{compatible with the $G$-action,}
\]
so it suffices to show that for every irreducible $\bQ_\ell$-representation $\tau$ of $G$, the multiplicities of $\tau$ and $\tau^*$ in $\cX_\ell(\Res_{F/K} (A_F))$ are equal. For this, we now suitably modify the proof of \cite{DD09c}*{Thm.~2.3}. Letting $d_\tau$ denote the dimension of any irreducible constituent of $\tau \tensor_{\bQ_\ell} \ov{\bQ}_\ell$, we consider
\[
\tst P_\tau \ce d_\tau \cdot \sum_{g \in G} \Tr(\tau(g))g \in \bZ_\ell[G].
\]
By the contravariance of $\cX_\ell(-)$ and by \cite{Ser77}*{\S2.6 Thm.~8 (ii) and \S12.2}, 
%Loc. cit. is written over a splitting field but that is no problem for us because $P_\tau$ is the sum of similar operators for absolutely irreducible constituents of $\tau$.
the operator $P_\tau$ kills every irreducible $G$-constituent $\tau' \not\simeq \tau$ of $\cX_\ell(\Res_{F/K} (A_F))$ and acts as scaling by $m_\tau \cdot \#G$ on every copy of $\tau$, where $m_\tau \in \bZ_{> 0}$ denotes the Schur index of $\tau$. Therefore, there exists a small $\ell$-adic neighborhood 
\[
U_\tau \quad  \text{of}  \quad m_\tau\cdot \#G \cdot 1_G + (p - 1)P_\tau \quad \text{in} \quad \bZ_\ell[G]
\]
such that for any $G$-stable $\bZ_\ell$-lattice $\Lambda_{\tau'}$ in an irreducible $G$-subrepresentation $\tau' \subset \cX_\ell(\Res_{F/K}(A_F))$ with $\tau' \not\simeq \tau$ (resp.,~with $\tau' \simeq \tau$), any $x \in U_\tau$ acts as the $m_\tau\cdot \#G$-multiple (resp.,~as the $\ell \cdot m_\tau \cdot\#G$-multiple) of a $\bZ_\ell$-automorphism of $\Lambda_{\tau'}$. Since, moreover, the $\bZ_\ell$-linear automorphism $\iota$ of $\bZ_\ell[G]$ determined by $\iota(g) = g\i$ carries $P_\tau$ to $P_{\tau^*}$ and is continuous, we may find an element
\[
\tst \Phi_\tau = \sum_{g \in G} x_{\tau,\, g} \cdot g \in \bZ[G] \qq \text{such that} \qq \Phi_\tau \in U_\tau \quad \text{and} \quad \iota(\Phi_\tau) \in U_{\tau^*}.
\]
Thanks to the $\ell \neq \Char K$ assumption we may assume, in addition, that the determinant of the $\bZ$-linear multiplication by $\Phi_\tau$ endomorphism of $\bZ[G]$ is prime to $\Char K$. Then the endomorphism $\phi_\tau$ of $\Res_{F/K}(A_F)$ determined by $\Phi_\tau$ is an isogeny of degree prime to $\Char K$, so \Cref{Mil-lem} gives the equality
\[
Q(\phi_\tau) = Q(\phi^\vee_\tau).
\]
Under a principal polarization, $\phi_\tau^\vee$ is identified with the endomorphism of $\Res_{F/K} (A_F)$ determined by $\iota(\Phi_\tau)$. Therefore, letting $\mathrm{mult}(\tau, V)$ denote the multiplicity of the irreducible $\bQ_\ell$-representation $\tau$ of $G$ in a $\bQ_\ell$-representation $V$ of $G$ (and likewise for $\tau^*$), we conclude from the construction of $\Phi_\tau$  that
\[\ba
\ord_\ell(Q(\phi_\tau)) &= \ord_\ell(m_\tau \cdot \#G) \cdot  \dim_{\bQ_\ell} \cX_\ell(\Res_{F/K} (A_F)) + \mathrm{mult}(\tau, \cX_\ell(\Res_{F/K} (A_F))) \cdot \dim_{\bQ_\ell} \tau, \\
\ord_\ell(Q(\phi_\tau^\vee)) &= \ord_\ell(m_{\tau^*} \cdot \#G) \cdot  \dim_{\bQ_\ell} \cX_\ell(\Res_{F/K} (A_F)) + \mathrm{mult}(\tau^*, \cX_\ell(\Res_{F/K} (A_F)))\cdot \dim_{\bQ_\ell} \tau^*.
\ea\]
Since $m_\tau = m_{\tau^*}$ and $\dim_{\bQ_\ell} \tau = \dim_{\bQ_\ell} \tau^*$, we obtain the desired equality
\[
\mathrm{mult}(\tau, \cX_\ell(\Res_{F/K} (A_F))) = \mathrm{mult}(\tau^*, \cX_\ell(\Res_{F/K} (A_F))). \qedhere
\]
\epf

\brem
It is desirable to remove the assumption $\ell \neq \Char K$ in \Cref{DD}. For this, the crux of the matter is to remove the degree restriction in \Cref{Mil-lem}.
\erem

With \Cref{DD} in hand, we proceed to \Cref{odd-gal}, whose  proof will use the following lemma:

\blem \lab{sel-gal}
For a finite Galois extension $F/K$ of global fields, a prime $\ell$, and an abelian variety $A$ over $K$, the map $\cX_\ell(A_{F}) \ra \cX_\ell(A)$ induced by restriction to $F$ supplies the second isomorphism in
\[
\cX_\ell(A_{F})^{G} \cong \cX_\ell(A_{F})_{G} \cong \cX_\ell(A), \qq \text{where $G \ce \Gal(F/K)$.}
\]
\elem

\bpf
The proof for elliptic curves and number fields, \cite{DD10}*{proof of Lemma 4.14}, extends: the spectral sequence 
\[
H^i(G, H^j_\fppf(F, A[\ell\I])) \Ra H^{i + j}_\fppf(K, A[\ell\I])
\]
 shows that $\#G$ kills the kernel and the cokernel of the map
 \[
 H^1_\fppf(K, A[\ell\I]) \ra H^1_\fppf(F, A[\ell\I])^G.
 \]
 Thus, $\#G$ also kills the kernel of the map
 \be \lab{map-map}
 \Sel_{\ell\I} A \ra(\Sel_{\ell\I} A_{F})^G.
 \ee
Moreover, $\#G$ kills $\Ker(H^1(K_v, A) \ra H^1(F_{w}, A))$ for all places $w$ of $F$ that extend a place $v$ of $K$, so $(\#G)^2$ kills the cokernel of the map \eqref{map-map}. In order to obtain the claim, it remains to pass to Pontryagin duals and to invert $\ell$.
\epf

\bcor \lab{odd-gal}
For an odd degree Galois extension $F/K$ of global fields, an abelian variety $A$ over $K$, and a prime $\ell$ different from $\Char K$, one has
\benum
\item $ \rk_\ell A \equiv \rk_\ell A_{F} \bmod 2$ and

\item $w(A) = w(A_{F})$.
\eenum
In particular, the $\ell$-parity conjecture holds for $A$ if and only if it holds for $A_{F}$.
\ecor

\bpf \hfill
\benum
\item Combine the argument of \cite{DD09c}*{Cor.~2.5} with \Cref{DD} and \Cref{sel-gal}.

\item The number field proof \cite{DD09b}*{A.2 (3)} also works for global fields. \qedhere  
\eenum
\epf

%\brem \lab{DD-app}
%Thanks to \Cref{DD} and \Cref{sel-gal}, all the results of \cite{DD09b}*{Appendix A} with the possible exception of \cite{DD09b}*{A.2 (5)} extend to global fields under the $p \neq \Char K$ restriction inherited from \Cref{DD}; the same applies to \cite{DD09c}*{Cor.~2.7}. We have excluded \cite{DD09b}*{A.2~(5)} due to its reliance on \cite{Roh11a}, which is written in characteristic $0$ context.
%\erem

\section{Replacing $\ell\I$-Selmer groups by $\ell$-Selmer groups} \lab{sel-sel}

To facilitate the Grothendieck--Ogg--Shafarevich input to our proof of \Cref{main}, in \Cref{rk-sel} we (implicitly) reformulate the $\ell$-parity conjecture by relating the $\ell\I$-Selmer rank and the $\ell$-Selmer rank. A suitable polarization is handy for this---without it, controlling the parity of $\dim_{\bF_\ell} (\Sha(A)_\nd[\ell])$ would become a major concern. In fact, even with it, this parity may vary in the $\ell = 2$ case, as Poonen and Stoll explain in \cite{PS99}. As far as the proof of \Cref{main} is concerned, the goal of \Cref{PS} and \Cref{CT-BC} below is to overcome this difficulty by proving that the said parity is even over every quadratic extension. \Cref{PS} is a slight improvement to the main result of \emph{op.~cit.}---without this improvement, in the $\ell = 2$ case of \Cref{main} we would be forced to restrict to principally polarized abelian varieties, which were the main focus of Poonen and Stoll.

\bpp[The Cassels--Tate pairing] \lab{CT}
For an abelian variety $B$ over a global field $F$, let 
\[
\langle\ ,\,\rangle\colon \Sha(B) \times \Sha(B^\vee) \ra \bQ/\bZ
\]
be the Cassels--Tate bilinear pairing. 
For a self-dual homomorphism $\lambda\colon B \ra B^\vee$, the pairing 
\[
\langle a, b\rangle_{\gL} \ce \langle a,\lambda(b)\rangle\quad\text{  for} \quad a, b \in \Sha(B) 
\]
is antisymmetric \cite{PS99}*{\S6, Cor.~6}. Therefore, if $\lambda$ is in addition an isogeny, then, for every prime $p \nmid 2\deg \lambda$, the pairing induced by $\langle\ , \, \rangle_{\gL}$ on the abelian group $\Sha(B)[p\I]$ is alternating. In this case, since $\In{\ ,\, }$ is nondegenerate modulo the divisible subgroups,\footnote{The first published complete proof of the fact that the left and the right kernels of $\In{\ ,\, }$ are the maximal divisible subgroups of $\Sha(B)$ and $\Sha(B^\vee)$, respectively, seems to be the combination of \cite{HS05}*{Thm.~0.2}, \cite{HS05e}, and \cite{GA09}*{Thm.~1.2} (although \emph{op.~cit.}~does not treat the prime to the characteristic torsion subgroups, its MathSciNet review remarks that the corresponding claim for such subgroups follows from the proof of \cite{HS05}*{Thm.~0.2}).} for every prime $p\nmid 2\deg \gL$ we have
\[
\dim_{\bF_p}(\Sha(B)_\nd[p]) \equiv 0 \bmod 2.
\]
For a self-dual homomorphism $\lambda\colon B \ra B^\vee$, as in \cite{PS99}*{\S4, Cor.~2}, we let
\[
c_{\lambda} \in \Sha(B^\vee)[2] \subset H^1(k, B^\vee)[2] \qq \text{be the image of} \qq \lambda \in (\NS B)(F).
\] 
In addition, for a self-dual isogeny  $\gL\colon B \ra B^\vee$  of odd degree, we let 
\[
\Sha(\gL)[2]\colon \Sha(B)[2] \isomto \Sha(B^\vee)[2]
\]
be the induced isomorphism and set
\[
c \ce (\Sha(\lambda)[2])\i( c_{\lambda}) \in \Sha(B)[2].
\]
For such $\gL$, 
\[
\tst \dim_{\bF_2}(\Sha(B)_\nd[2]) \mod 2 \qq \text{is governed by} \qq \In{c, c}_{\lambda} \in \{ 0, \f{1}{2}\};
\]
indeed,  Poonen and Stoll observed this in \cite{PS99}*{\S6,~Thm.~8} in the case when $\lambda$ is a principal polarization, and the general case follows from their argument:
\epp

\bthm \lab{PS}
Let $\lambda\colon B \ra B^\vee$ be a self-dual isogeny  of odd degree $d$.
\benum
\item \lab{PS-a}
If $\In{c, c}_{\gL} = 0$, then $\#\Sha(B)_\nd[n]$ is a perfect square for every $n \in \bZ_{\ge 1}$ prime to $d$.

\item
If $\In{c, c}_{\gL} = \f{1}{2}$, then $\#\Sha(B)_\nd[n]$ is twice a perfect square for every even $n \in \bZ_{\ge 1}$ prime to $d$ and is a perfect square for every odd $n \in \bZ_{\ge 1}$ prime to $d$.
\eenum
\ethm

\bpf
The proof of \cite{PS99}*{\S6,~Thm.~8} given for the case $d = 1$ continues to work if throughout that proof one replaces $\Sha_{\nd}$ by its subgroup $\Sha_{\nd}(d\pr)$ consisting of the elements of order prime to $d$.
\epf

\brem
Under the assumptions of \Cref{PS}, the replacement indicated in the proof also gives the following extension of  \cite{PS99}*{\S6, Cor.~9}: if $\Sha(B)_\nd$ is finite, then there exists a finite abelian group $T$ such that, letting $(-)(d')$ denote the ``prime to $d$'' subgroup,
\[\ba
\text{in the case}\qq  \In{c, c}_{\gL} = 0,\, \qq &\x{one has} \qq \Sha(B)_\nd(d') \simeq T \times T; \\ 
 \x{in the case} \quad \ \ \, \In{c, c}_{\gL} = \f{1}{2}, \qq &\x{one has} \qq \Sha(B)_\nd(d') \simeq \bZ/2\bZ \times T \times T.
\ea\]
\erem

Due to \Cref{CT-BC} below, the obstruction $\In{c, c}_\lambda$ vanishes over an even degree extension. \Cref{CT-BC} was pointed out to us by Bjorn Poonen and was also observed by Adam Morgan.

\blem \lab{CT-BC}
For a finite extension $F\pr/F$ of global fields, there is a commutative diagram
\[
\xymatrix{
\Sha(B) \times \Sha (B^\vee) \ar[d]^-{\Res} \ar[r]^-{\langle\, ,\, \rangle} & \bQ/\bZ \ar[d]^-{[F\pr : F]} \\
\Sha(B_{F\pr}) \times \Sha (B_{F\pr}^\vee) \ar[r]^-{\langle\, ,\, \rangle} & \bQ/\bZ.
}
\]
\elem

\bpf
Combine the definition \cite{PS99}*{\S3.1} of the pairings with the well-known commutativity of the diagram
\[
\xymatrix{
\Br(F_v) \ar[d]_{\Res}\ar[r]^-{\inv_v} & \bQ/\bZ \ar[d]^-{[F_{v\pr}\pr : F_v]} \\
\Br(F\pr_{v\pr})  \ar[r]^-{\inv_{v\pr}} & \bQ/\bZ
}
\]
for a finite extension $F\pr_{v\pr}/F_v$ of local fields (for the commutativity, see \cite{Ser67}*{\S1.1, Thm.~3}).
\epf

We are ready to reduce to working with $\ell$-Selmer groups instead of $\ell\I$-Selmer groups in \Cref{main}.

\bcor \lab{rk-sel}
For a prime $\ell$, if an abelian variety $B$ over a global field $F$ has a polarization of degree prime to $\ell$, then for every quadratic extension $F'/F$ one has
\[
\rk_\ell B_{F'} \equiv \dim_{\bF_\ell} \Sel_\ell B_{F'} - \dim_{\bF_\ell} B[\ell](F') \bmod 2.
\]
\ecor

\bpf
For any prime $\ell$ and any abelian variety $A$ over a global field $K$ one has
\[
\rk_\ell A = \rk A + \dim_{\bF_\ell} \Sha(A)[\ell] - \dim_{\bF_\ell} \Sha(A)_\nd[\ell] = \dim_{\bF_\ell} \Sel_\ell A - \dim_{\bF_\ell} A[\ell](K) - \dim_{\bF_\ell} \Sha(A)_\nd[\ell].
\]
Moreover, for $A = B_{F'}$, the dimension $\dim_{\bF_\ell} \Sha(B_{F'})_\nd[\ell]$ is even by \Cref{PS}~\ref{PS-a} and \Cref{CT-BC}. The desired congruence follows.
\epf

\brem
Even though the proof continues to work, no generality is gained by requiring a self-dual isogeny instead of a polarization in \Cref{rk-sel} (or \Cref{main}): for an $n \in \bZ_{> 0}$, if an abelian variety $B$ over a field $F$ has a self-dual isogeny  of degree prime to $n$, then it also has a polarization of degree prime to $n$. Indeed, one knows that the degree function\footnote{The degree of a self-dual homomorphism that is not an isogeny is defined to be $0$.} $\deg\colon (\NS B)(F) \ra \bZ$ is a polynomial with rational coefficients on the lattice $(\NS B)(F)$; 
% Compare Grieve "Index conditions ..." section 4.1 and/or Moonen, van der Geer "Abelian varieties" discussion before 9.11
consequently, $\deg$ modulo $n$ is translation invariant with respect to a sublattice, and it remains to note that the cone of polarizations spans $(\NS B)(F)$.
% For instance, by Hartshorne, Ex. II.7.5 (b)
\erem

\section{Local root numbers in unramified extensions}\lab{app}

The goal of this section is \Cref{st-root}, which details the behavior of the local root number of an abelian variety upon an unramified extension of degree $m$ of the nonarchimedean local base field. In fact, this behavior manifests itself for a wider class of representations than those coming from abelian varieties, as we observe in \Cref{st-WD}. To summarize, for every representation in this class the local root number ``stabilizes'' upon unramified base change of sufficiently divisible degree to a value determined by the parity of the conductor.\footnote{We do not use conductor ideals, so `conductor' abbreviates what some authors call `conductor exponent.'} Such behavior, which is crucial for our proof of \Cref{main}, seems not to have been pointed out previously.

Throughout \S\ref{app}, we let $k$ be a nonarchimedean local field and let $\fo$, $\bF$, $p$, and $k^s$ be its ring of integers, its residue field, its residue characteristic, and a choice of a separable closure, respectively. We denote the unramified subextension of $k^s/k$ of degree $m$ and its ring of integers by $k_m$ and $\fo_m$, respectively. We let $W(k^s/k)$ and $I$ denote the Weil group and its inertia subgroup. We denote a geometric Frobenius in $W(k^s/k)$ by $\Frob_k$. For a field $F$ with $\Char F \neq p$, we let 
\[
\abs{\cdot}_k\colon W(k^s/k) \ra F^\times
\]
be the unramified character characterized by the equality $\abs{\Frob_k}_k = (\#\bF)\i$; for an integer $n$ and a representation $V$ of $W(k^s/k)$ over $F$, we set $V(n) \ce V \tensor_F \abs{\cdot}_k^n$, where the second factor denotes a copy of $F$ on which $W(k^s/k)$ acts through the $n\th$ power of $\abs{\cdot}_k$.

\bpp[$\eps$-factors of Weil--Deligne representations] \lab{eps-WD}
For a field $F$ with $\Char F \neq p$, a \emph{Weil--Deligne representation} of $W(k^s/k)$ over $F$ is a pair 
\[
\rho\pr = (\rho, N)
\]
that consists of
\begin{itemize}
\item a finite dimensional representation $\rho$ of $W(k^s/k)$ over $F$ such that the restriction of $\rho$ to some open subgroup of $I$ is trivial, and

\item a $W(k^s/k)$-homomorphism $N\colon \rho \ra \rho(-1)$. 
\end{itemize}
Subject to the choices of a nontrivial additive character $\psi\colon k \ra F^\times$ and a nonzero $F$-valued Haar measure $dx$ on $(k, +)$, the $\eps$-factor of $\rho\pr$ is defined by
\be\lab{eps-WD-def}
\eps(\rho\pr, \psi, dx) \ce \eps_0(\rho, \psi, dx) \det(-\Frob_k\, |\, (\Ker N)^{I})\i,
\ee
where for the appearing $\eps_0$-factor as well as for the definitions of an additive character and an $F$-valued Haar measure we refer to \cite{Del73}*{\S6} (or to \cite{Ces16a}*{1.1 and \S\S2.3--4}). The \emph{Artin conductor} of $\rho\pr$ is
\be\lab{art-cond}
a(\rho\pr) \ce \Sw \rho + \dim_F \rho - \dim_F (\Ker N)^{I};
\ee
for the definition of the Swan conductor $\Sw \rho$, see \cite{Del73}*{\S6.2} or \cite{Ces16a}*{\S2.9}.
\epp

Before proceeding, for later use we record the following lemma about conductors of abelian varieties.

\blem \lab{cond-comp}
Let $B \ra \Spec k$ be an abelian variety, let $a(B)$ be its conductor exponent, let $\cB_{\bF}$ be the special fiber of the N\'{e}ron model of $B$, and let $\Phi$ be the component group scheme of $\cB_{\bF}$. For every prime $\ell$ different from $\Char \bF$, one has
\[
a(B) = a(B[\ell]) + \dim_{\bF_\ell}(\Phi[\ell](\ov{\bF})),
\]
where the Artin conductor $a(B[\ell])$ is defined by \eqref{art-cond} {\upshape(}with $N = 0${\upshape)}.
\elem

\bpf
Let $(V_\ell B)^{\mathrm{ss}}$ be the semisimplification of the $\ell$-adic Tate module of $B$. Then 
\[\ba
a(B) &= \Sw((V_\ell B)^{\ss}) + \dim_{\bQ_\ell}(V_\ell B) - \dim_{\bQ_\ell} (V_\ell B)^{I}, \\
a(B[\ell]) &= \Sw(B[\ell]) + \dim_{\bF_\ell}  (B[\ell]) - \dim_{\bF_\ell}(B[\ell]^{I}).
\ea\]
The identification of $B[\ell]^{I}$ and $(V_\ell B)^{I}$ with $\cB_{\bF}[\ell]$ and $V_\ell (\cB_{\bF})$ explained in \cite{ST68}*{Lemma~2} gives
\[
\dim_{\bQ_\ell} (V_\ell B)^{I} = \dim_{\bQ_\ell} V_\ell(\cB_\bF) = \dim_{\bF_\ell}( \cB_\bF[\ell]) - \dim_{\bF_\ell} (\Phi[\ell](\ov{\bF})) = \dim_{\bF_\ell} (B[\ell]^{I}) - \dim_{\bF_\ell} (\Phi[\ell](\ov{\bF})),
\]
so it remains to note that $\Sw((V_\ell B)^{\ss}) = \Sw(B[\ell])$ because the Swan conductor is additive and is compatible with reduction mod $\ell$.
\epf

Returning to the setup of \S\ref{eps-WD}, we turn to the analysis of the epsilon factor of $\rho'|_{k_m}$.

\bprop \lab{eps-gr}
In the setup of {\upshape\S\ref{eps-WD}}, for the restriction $\rho\pr|_{k_m}$ of $\rho\pr$ to $W(k^s/k_m)$ one has
\[
\eps(\rho\pr|_{k_m}, \psi \circ \Tr_{k_m/k}, dx_{m}) = \begin{cases} \eps(\rho\pr, \psi, dx)^m,\quad\quad\quad\quad\, \text{if $m$ is odd,}  \\ (-1)^{a(\rho\pr)} \eps(\rho\pr, \psi, dx)^m, \, \text{ if $m$ is even.} \end{cases}
\]
Here $dx_m$ denotes the Haar measure on $(k_m, +)$ characterized by $\int_{\fo_{m}} dx_{m} = (\int_{\fo} dx)^m$. 
\eprop

\bpf
For the $\eps_0$-factor appearing in \eqref{eps-WD-def}, the inductivity in degree $0$ gives the equality
\be\lab{A}
\eps_0(\rho|_{k_m}, \psi \circ \Tr_{k_m/k}, dx_m) = \eps_0(\b{1}_{k_m}, \psi \circ \Tr_{k_m/k}, dx_m)^{\dim_F \rho} \cdot \f{\eps_0((\Ind_{k_m}^k \b{1}_{k_m}) \tensor \rho, \psi, dx)}{\eps_0(\Ind_{k_m}^k \b{1}_{k_m} , \psi, dx)^{\dim_F \rho}}.
\ee
Since $\Ind_{k_m}^k \b{1}_{k_m}$ is unramified, \cite{Del73}*{5.5.3} (or \cite{Ces16a}*{3.2.2} for general $F$) simplifies the fraction~to
\be\lab{B}
\det(\Ind_{k_m}^k \b{1}_{k_m})(\Frob_k)^{\Sw \rho} \cdot \f{\eps_0(\rho, \psi, dx)^m}{\eps_0(\b{1}_{k} , \psi, dx)^{m\dim_F \rho}} = (-1)^{(m - 1) \Sw \rho} \cdot \f{\eps_0(\rho, \psi, dx)^m}{\eps_0(\b{1}_{k} , \psi, dx)^{m\dim_F \rho}}.
\ee
Let $n(\psi)$ denote the largest integer $n$ such that $\psi|_{\pi^{-n}\fo} = 1$, where $\pi \in \fo$ is a uniformizer. Since $k_m/k$ is unramified, $n(\psi \circ \Tr_{k_m/k}) = n(\psi)$ by \cite{Del73}*{\S4.11}; we use this in the following computation:
\be\lab{C}
\eps_0(\b{1}_{k} , \psi, dx)^{m} = \p{-(\#\bF)^{n(\psi)} \cdot \int_{\fo} dx}^m = (-1)^{m - 1}\eps_0(\b{1}_{k_m}, \psi \circ \Tr_{k_m/k}, dx_m).
\ee
The equations \eqref{A}, \eqref{B}, and \eqref{C} combine to give the equality
\be\lab{D}
\eps_0(\rho|_{k_m}, \psi \circ \Tr_{k_m/k}, dx_m) = (-1)^{(m - 1)(\Sw \rho + \dim_F \rho)} \eps_0(\rho, \psi, dx)^m.
\ee
It remains to put \eqref{D} together with the evident equality
\[
\det(-\Frob_{k_m}\, |\, (\Ker N)^{I})\i = (-1)^{-(m - 1)\dim_F (\Ker N)^{I}} \det(-\Frob_{k}\, |\, (\Ker N)^{I})^{-m}. \qedhere
\]
\epf

\bpp[Root numbers] \lab{root}
Assume that $F = \bC$ and $\int_{\fo} dx \in \bR^+$ in \S\ref{eps-WD}. The \emph{root number} of $\rho\pr$ is
\[
w(\rho\pr, \psi) \ce \f{\eps(\rho\pr, \psi, dx)}{\abs{\eps(\rho\pr, \psi, dx)}}.
\]
It does not depend on the choice of $dx$ as long as $\int_{\fo} dx \in \bR^+$. If $\det \rho$ is $\bR^+$-valued, then $w(\rho\pr, \psi)$ does not depend on the choice of $\psi$ either, thanks to the formula \cite{Del73}*{5.4}. In this case we abbreviate $w(\rho\pr, \psi)$ by $w(\rho\pr)$. If $B$ is an abelian variety over $k$, then
\[
\textstyle{\bigwedge^{2g}} H^1_\et(B, \bQ_\ell) \cong H^{2g}_\et(B, \bQ_\ell) \cong \bQ_\ell(-g),
\]
% The first iso. is induced by cup product (see Milne "Abelian varieties", 15.1 (b), so it is Galois equivariant. The second one is standard for smooth geom. connected varieties.
so the independence of $\psi$ is witnessed if $\rho\pr$ is the complex Weil--Deligne representation $\sigma_B\pr$ that one associates to $H^1_\et(B, \bQ_\ell) \cong (V_\ell B)^*$ for a prime $\ell$ different from $p$ using the Grothendieck quasi-unipotence theorem and an embedding $\iota\colon \bQ_\ell \hra \bC$. By \cite{Sab07}*{1.15}, the isomorphism class of $\sigma_B\pr$ does not depend on $\ell$ and $\iota$,\footnote{\emph{Loc.~cit.}~does not use its additional $\Char k = 0$ assumption in the proof. Also, we bypass this issue by analyzing the right hand side of \eqref{l-par-conj} through the second case of \Cref{st-root}, the proof of which works for every $\ell$ and $\iota$.} and hence neither does the root number of $B$ defined by 
\[
w(B) \ce w(\sigma_B\pr).
\]
Due to the Weil pairing, the presence of a polarization of $B$, and \cite{Del73}*{5.7.1}, one has $w(B) \in \{ \pm 1\}$.
\epp

\bcor\lab{st-WD}
For a Weil--Deligne representation $\rho\pr$ of $W(k^s/k)$ over $\bC$ such that $\det\rho$ is $\bR^+$-valued and $w(\rho\pr)$ is an $m_0\th$ root of unity,
\[
w(\rho\pr|_{k_m}) = (-1)^{a(\rho\pr)}
\]
for every even $m$ divisible by $m_0$.
\ecor

\bpf
This follows from \Cref{eps-gr}.
\epf

\bcor \lab{st-root}
Let $B$ be an abelian variety over $k$, and let $a(B)$ be its conductor exponent. Then
\[
w(B_{k_m}) = \begin{cases} w(B),\quad \quad\  \text{ if $m$ is odd,} \\  (-1)^{a(B)},\quad \text{ if $m$ is even.}  \end{cases}
\]
\ecor

\bpf
Combine \Cref{eps-gr} and the equality $a(B) = a(\sigma_B\pr)$ that results from the definitions (for which one may consult \cite{Ser70}*{\S2}).
\epf

\brem
For elliptic curves, excluding the troublesome additive reduction case if $p \le 3$, one may also prove \Cref{st-root} by the means of explicit case-by-case formulae for $w(B)$ and $a(B)$.
\erem

\section{Arithmetic duality generalities: comparing Selmer sizes modulo squares} \lab{ad-pre}

The main goal of this section is to prove \Cref{beyond-KMR}, which in \S\ref{ade} will specialize to the arithmetic duality input needed for the proof of \Cref{main}. \Cref{beyond-KMR} generalizes \cite{KMR13}*{Thm.~3.9} to the case of commutative self-dual finite group schemes over global fields from the case of self-dual $2$-dimensional $\bF_p$-vector space group schemes over number fields. Although its proof is loosely modeled on that of \emph{loc.~cit.},~modifications are necessary due to the possibility that $\Char F \mid \# \cG$, when various cohomology groups are no longer finite. The simpler case of \Cref{beyond-KMR} when $\Char F \nmid \#\cG$ suffices for the proof of \Cref{main}, but it seems unnatural to confine the general techniques in this way. Consequently, \Cref{dual} does not exclude the more subtle cases when $\Char F \mid n$.

In the buildup to \Cref{beyond-KMR} we follow an axiomatic approach by introducing further assumptions as we need them. This way, in \Cref{beyond-MR} we arrive at a generalization of \cite{MR07}*{Prop.~1.3 (i)} that removes the self-duality, $\bF_p$-vector space, and number field assumptions from \emph{loc.~cit.}

In this section and in \S\ref{ade} all the cohomology groups  are fppf. Identifications with \'{e}tale or Galois cohomology are implicit. Likewise implicit is the Tate modification: if $v$ is archimedean, then we write $H^i(K_v, -)$ for $\wh{H}^i(K_v, -)$ (this has no effect if $i \ge 1$).

\bpp[The basic setup] \lab{setup}
Let $F$ be a global field. If $\Char F = 0$, let $S$ be the spectrum of the ring of integers of $F$; if $\Char F > 0$, let $S$ be the connected smooth proper curve over a finite field such that the function field of $S$ is $F$. Let $U \subset S$ be a nonempty open subscheme. We denote by $v$ a place of $F$ and identify the nonarchimedean $v$ with the closed points of $S$; writing $v\not \in U$ signifies that $v$ does not correspond to a closed point of $U$ (and hence could be archimedean). 

Let 
\[
\cG \ra U \qq \text{and} \qq \cH \ra U
\]
be commutative finite flat group schemes, and suppose that there is a perfect bilinear pairing 
\[
\cG \times_U \cH \xra{b} \bG_m \qqq \text{that identifies $\cG$ and $\cH$ as Cartier duals.}
\]
The cohomology groups $H^1(U, \cG)$ and $H^1(U, \cH)$ are ``cut out by local conditions,'' i.e., as noted in \cite{Ces16c}*{4.3}, the squares
\[
\xymatrix{
H^1(U, \cG) \ar[d] \ar@{^(->}[r] & H^1(F, \cG)  \ar[d]    	& & H^1(U, \cH) \ar[d] \ar@{^(->}[r] & H^1(F, \cH)  \ar[d]    \\
\prod\limits_{v\in U} H^1(\cO_v, \cG) \ar@{^(->}[r] & \prod\limits_{v\in U} H^1(F_v, \cG), & & \prod\limits_{v\in U} H^1(\cO_v, \cH) \ar@{^(->}[r] & \prod\limits_{v\in U} H^1(F_v, \cH)
}
\]
are Cartesian. The main result of this section, \Cref{beyond-KMR}, investigates further subgroups cut out by also imposing local conditions at all $v\not\in U$. Its proof hinges on, among other things, the Tate--Shatz local duality \cite{Sha64}*{Duality theorem on p.~411} (alternatively, \cite{Mil06}*{I.2.3,~I.2.13~(a),~III.6.10}), which says that for every place $v$ and integer $i$ the cup product pairing
\be\lab{Shatz}
\xymatrix{
H^i(F_v, \cG) \times H^{2 - i}(F_v, \cH) \ar[r] &H^2(F_v, \bG_m) \ar@{^(->}[r]^-{\inv_v} &\bQ/\bZ
}
\ee
that uses $b$ identifies $H^i(F_v, \cG)$ and $H^{2-i}(F_v, \cH)$ as Pontryagin duals of each other. The Pontryagin duality in question is that of locally compact Hausdorff abelian topological groups---see \cite{Ces15b}*{3.1--3.2 and 1.4} for the definition and properties of the topology on these cohomology groups. These groups are finite and discrete if $\Char F_v \nmid \#\cG$; this case suffices for the proof of \Cref{main}.
\epp

\beg \lab{main-eg}
Our main case of interest in the setup of \S\ref{setup} is when $B \ra \Spec F$ and $B^\vee \ra \Spec F$ are dual abelian varieties, $\cB \ra S$ and $\cB^\vee \ra S$ are their N\'{e}ron models, $U \subset S$ is such that $\cB_U \ra U$ and $\cB_U^\vee \ra U$ are abelian schemes, and $\cG = \cB[n]_U$, $\cH = \cB^\vee[n]_U$ for some $n \in \bZ_{> 0}$. Cartier--Nishi duality \cite{Oda69}*{Thm.~1.1} supplies the pairing $b$ in this case.
\eeg

The following lemma is crucial for the arithmetic duality results derived below.

\blem \lab{im-orth}
In the setup of {\upshape\S\ref{setup}}, for every integer $i$ the images of the pullback maps
\[
 H^i(U, \cG) \xra{\loc^i(\cG)} \textstyle\bigoplus_{v\not\in U} H^i(F_v, \cG) \qq \text{and} \qq H^{2 - i}(U, \cH) \xra{\loc^{2 - i}(\cH)} \bigoplus_{v \not\in U} H^{2 - i}(F_v, \cH)
\]
are orthogonal complements under the sum of the pairings \eqref{Shatz}.
\elem

\bpf
In the case when $\#\cG \in \Gamma(U, \cO_U^\times)$, the Poitou--Tate sequence gives the claim once one explicates its morphisms and interprets the global cohomology groups as Galois cohomology with restricted ramification (for this interpretation consult, e.g., \cite{Mil06}*{II.2.9}). To treat the general case we will use an extension of the Poitou--Tate sequence, namely, the compactly supported flat cohomology~sequence
\be \lab{ccch}
\dotsb \ra H^i_c(U, \cG) \ra H^i(U, \cG) \xra{\loc^i(\cG)} \textstyle\bigoplus_{v \not\in U} H^i(F_v, \cG) \xra{\delta_{c}^i(\cG)} H^{i + 1}_c(U, \cG) \ra \dotsb
\ee
of \cite{Mil06}*{III.0.4 (a)}. The pairings in the diagram 
\be\tag{$\dagger$}\ba\lab{pair}
\xymatrix{
H^i(U, \cG) \ar[d]_-{\loc^i(\cG)} \ar@{}[r]|-{\bigtimes} & H^{3- i}_c(U, \cH) \ar[rrr]^{\text{\cite{Mil06}*{III.3.2 and III.8.2}}} &&& H^3_c(U, \bG_m) \ar[r]^-{\tr} &\bQ/\bZ \ar@{=}[d] \\
\bigoplus_{v\not\in U} H^i(F_v, \cG) \ar@{}[r]|-{\bigtimes} & \bigoplus_{v\not\in U} H^{2 - i}(F_v, \cH) \ar[u]^-{\delta^{2 - i}_c(\cH)} \ar[rrr]^{\sum_v \eqref{Shatz}} &&& \bigoplus_{v\not\in U} H^2(F_v, \bG_m) \ar[u]^-{\delta^2_c(\bG_m)} \ar[r]^-{\sum_v \inv_v} & \bQ/\bZ
}
\ea\ee
are perfect, so the exactness of \eqref{ccch} reduces the desired claim to the commutativity of both squares of \eqref{pair} (and of their analogues with $\cG$ and $\cH$ interchanged and $i$ replaced by $2 - i$). The right square of \eqref{pair} commutes by \cite{Mil06}*{(b) in the beginning of II.\S3}.

The cup product \eqref{Shatz} in the bottom left arrow of \eqref{pair} agrees with the $\Ext$-product as in \cite{GH71}*{3.1}, which in turn agrees with the Yoneda edge product as in \cite{GH70}*{4.5}, so the left square of \eqref{pair} is identified with the square
\be\tag{\ddag}\ba\lab{pair-pair}
\xymatrix{
\Ext^i_U(\cH, \bG_m) \ar[d] \ar@{}[r]|-{\bigtimes} & H^{3- i}_c(U, \cH) \ar[rr]^{\text{\cite{Mil06}*{III.0.4 (e)}}} && H^3_c(U, \bG_m)  \\
\bigoplus_{v\not\in U} \Ext^i_{F_v}(\cH, \bG_m) \ar@{}[r]|-{\bigtimes} & \bigoplus_{v\not\in U} H^{2 - i}(F_v, \cH) \ar[u]^-{\delta^{2 - i}_c(\cH)} \ar[rr] && \bigoplus_{v\not\in U} H^2(F_v, \bG_m) \ar[u]^-{\delta^2_c(\bG_m)}.
}
\ea\ee
The desired commutativity of \eqref{pair-pair} then results from the definitions and from the inspection of the proof of \cite{Mil06}*{III.0.4~(e)}: if one fixes injective resolutions 
\[
\cH \ra I^\bullet(\cH) \qq \text{and} \qq \bG_m \ra I^\bullet(\bG_m)
\]
over $U$ and interprets elements of $\Ext^i(\cH, \bG_m)$ and $\bigoplus_{v\not\in U} H^{2 - i}(F_v, \cH)$ as homotopy classes of maps 
\[
\tst I^\bullet(\cH) \xra{a} I^\bullet(\bG_m)[i] \qq \text{and} \qq \bZ \xra{d} \Gamma\p{\bigoplus_{v\not\in U} F_v, I^\bullet(\cH)|_{\bigoplus_{v\not \in U} F_v}}[2- i],
\]
then both ways to pair $a$ and $d$ in \eqref{pair-pair} result in the element of $H^3_c(U, \bG_m)$ that is represented by the homotopy class of the map
\[
\xymatrix{
\bZ \ar[rrrrr]^-{\p{0,\, \Gamma(\bigoplus_{v\not\in U} F_v,\, a|_{\bigoplus_{v\not\in U} F_v})[2 - i]\, \circ\, d}} &&&&&\textstyle \Gamma\p{U, I^\bullet(\bG_m)}[3] \oplus \Gamma\p{\bigoplus_{v\not\in U} F_v, I^\bullet(\bG_m)|_{\bigoplus_{v\not \in U} F_v}}[2].
} \qedhere
\]
% I leave it to Milne to unravel whatever is happening for the archimedean places (including what the actual definitions are); another minor problem is that in the number field case he uses fraction fields of the Henselizations, whereas for function fields he uses completions; annoying!
\epf

\bpp[Local conditions at $v\not\in U$] \lab{not-in-U}
In the setup of \S\ref{setup}, suppose that for every $j \in \{1, 2\}$ and $v \not \in U$ we have open compact subgroups 
\[
\Sel^j(\cG_{F_v}) \subset H^1(F_v, \cG) \qq \text{and} \qq \Sel^j(\cH_{F_v}) \subset H^1(F_v, \cH)
\]
that are orthogonal complements under \eqref{Shatz}. For $j \in \{1, 2\}$, define the Selmer groups $\Sel^j(\cG)$ and $\Sel^j(\cH)$ by requiring the sequences
\[\ba
0 \ra &\Sel^j(\cG) \ra H^1(U, \cG) \ra \textstyle\bigoplus_{v\not\in U} H^1(F_v, \cG)/\Sel^j(\cG_{F_v}),\\
0 \ra &\Sel^j(\cH) \ra H^1(U, \cH) \ra \textstyle\bigoplus_{v\not\in U} H^1(F_v, \cH)/\Sel^j(\cH_{F_v})
\ea\]
to be exact. For $v\not\in U$, \cite{HR79}*{24.10} gives further orthogonal complements 
\[
\Sel^{1 + 2}(\cG_{F_v}) \ce \Sel^1(\cG_{F_v}) + \Sel^2(\cG_{F_v})\qq \text{and} \qq \Sel^{1 \cap 2}(\cH_{F_v}) \ce \Sel^1(\cH_{F_v}) \cap \Sel^2(\cH_{F_v}),
\]
and one defines the corresponding Selmer groups $\Sel^{1 + 2}(\cG)$ and $\Sel^{1 \cap 2}(\cH)$ by the exactness of the sequences
\[\ba
0 \ra &\Sel^{1 + 2}(\cG) \ra H^1(U, \cG) \ra \textstyle\bigoplus_{v\not\in U} H^1(F_v, \cG)/\Sel^{1 + 2}(\cG_{F_v}),\\
0 \ra &\Sel^{1 \cap 2}(\cH) \ra H^1(U, \cH) \ra \textstyle\bigoplus_{v\not\in U} H^1(F_v, \cH)/\Sel^{1 \cap 2}(\cH_{F_v}).
\ea\]
\epp

\beg \lab{eg-loc-c}
In \Cref{main-eg} the local conditions of most interest to us arise when one takes the images of the local Kummer maps for $\Sel^1$ and, under appropriate restrictions, their fppf cohomological counterparts $H^1(\cO_v, \cB[n])$ and $H^1(\cO_v, \cB^\vee[n])$ for $\Sel^2$. In \S\ref{ade}, we will justify that the results of \S\ref{ad-pre} can be applied in this setting to compare $\#\Sel_n B$ and $\#H^1(S, \cB[n])$ modulo squares.
\eeg

\bprop \lab{beyond-MR}
In the setup of {\upshape\S\ref{not-in-U}}, the Selmer groups $\Sel^j(\cG)$ and $\Sel^j(\cH)$ for $j \in \{1, 2\}$, as well as $\Sel^{1 + 2}(\cG)$ and $\Sel^{1\cap 2}(\cH)$, are finite and 
\[
\#\p{\f{\Sel^{1 + 2}(\cG)}{\Sel^1(\cG)}} \cdot \#\p{\f{\Sel^1(\cH)}{\Sel^{1 \cap 2}(\cH)}} = \prod_{v\not\in U} \#\p{\f{\Sel^{1 + 2}(\cG_{F_v})}{\Sel^1(\cG_{F_v})}} = \prod_{v\not\in U} \#\p{\f{\Sel^1(\cH_{F_v})}{\Sel^{1\cap 2}(\cH_{F_v})}}.
\]
\eprop

\bpf
The finiteness is a special case of \cite{Ces16b}*{3.2} (whose proof uses Lemma 5.3 but no other results of this paper).

Due to \Cref{im-orth} and the choice of $\Sel^1(\cG_{F_v})$ and $\Sel^1(\cH_{F_v})$, \cite{HR79}*{24.10} shows that the subgroups
\[
\im(\loc^1(\cG))\,\textstyle{+} \bigoplus_{v\not\in U} \Sel^1(\cG_{F_v})  \subset \bigoplus_{v\not\in U} H^1(F_v, \cG)  \quad\text{and}\quad \im(\loc^1(\cH)|_{\Sel^1(\cH)}) \subset \bigoplus_{v\not\in U} H^1(F_v, \cH)
\]  
are orthogonal complements. Therefore, so are the subgroups
\[
\textstyle{} \f{H^1(U, \cG)}{\Sel^1(\cG)} \subset \bigoplus_{v\not\in U} \f{H^1(F_v, \cG)}{\Sel^1(\cG_{F_v})}  \quad\text{and}\quad \im(\loc^1(\cH)|_{\Sel^1(\cH)}) \subset \bigoplus_{v\not\in U} \Sel^1(\cH_{F_v}).
\]
Likewise, the subgroups
\[
\textstyle{} \bigoplus_{v\not\in U}\f{\Sel^{1+ 2}(\cG_{F_v})}{\Sel^1(\cG_{F_v})} \subset \bigoplus_{v\not\in U} \f{H^1(F_v, \cG)}{\Sel^1(\cG_{F_v})}  \quad\text{and}\quad \bigoplus_{v\not\in U} \Sel^{1\cap 2}(\cH_{F_v}) \subset \bigoplus_{v\not\in U} \Sel^1(\cH_{F_v})
\]
are also orthogonal complements. By combining the last two claims we deduce that the subgroups
\[
\textstyle{} \f{\Sel^{1 + 2}(\cG)}{\Sel^1(\cG)}  \subset \bigoplus_{v\not\in U} \f{H^1(F_v, \cG)}{\Sel^1(\cG_{F_v})}  \quad\text{and}\quad \im(\loc^1(\cH)|_{\Sel^1(\cH)}) + \bigoplus_{v\not\in U} \Sel^{1\cap 2}(\cH_{F_v}) \subset \bigoplus_{v\not\in U} \Sel^1(\cH_{F_v})
\]
are orthogonal complements, too, and hence so are the subgroups
\[
\textstyle{} \f{\Sel^{1 + 2}(\cG)}{\Sel^1(\cG)}  \subset \bigoplus_{v\not\in U} \f{\Sel^{1+ 2}(\cG_{F_v})}{\Sel^1(\cG_{F_v})}  \quad\text{and}\quad \f{\Sel^1(\cH)}{\Sel^{1\cap 2}(\cH)} \subset \bigoplus_{v\not\in U} \f{\Sel^1(\cH_{F_v})}{\Sel^{1 \cap 2}(\cH_{F_v})}.
\]
Both $\Sel^1(\cH_{F_v})$ and $\Sel^{1 \cap 2}(\cH_{F_v})$ are open and compact, so the groups in the last display are finite and the conclusion follows from the fact that dual finite abelian groups have the same cardinality.
\epf

\bpp[Breaking the symmetry] \lab{break}
In the setup of \S\ref{not-in-U}, suppose that one has the following data.
\begin{itemize}
\item
An isomorphism 
\[
\qq \theta\colon \cG \isomto \cH
\]
of $U$-group schemes such that $H^1(F_v, \theta)$  identifies the subgroup
\[
\qqq \Sel^j(\cG_{F_v}) \subset H^1(F_v, \cG) \quad \text{with} \quad \Sel^j(\cH_{F_v}) \subset H^1(F_v, \cH) \quad\quad \text{for every $v\not\in U$ and $j \in \{1, 2\}$,}
\]
so that, consequently, the isomorphism $H^1(U, \theta)$ identifies the subgroup
\[
\quad \Sel^j(\cG) \subset H^1(U, \cG)\quad  \text{with}\quad  \Sel^j(\cH) \subset H^1(U, \cH) \qq \text{for every $j \in \{1, 2\}$.}
\]

\item 
For each $v\not\in U$, a map $q_v \colon H^1(F_v, \cG) \ra \bQ/\bZ$ subject to the following requirements.
\benumr
\item \lab{quad-form}
One has $q_v(ax) = a^2 q_v(x)$ for all $a \in \bZ$ and $x \in H^1(F_v, \cG)$.

\item \lab{bil-pair}
The map 
\[
\qqq (x, y) \mapsto q_v(x + y) - q_v(x) - q_v(y)
\]
agrees with the bilinear form $\In{x, y}_v$ defined by the commutativity of the diagram
\be\ba\lab{pair-comp}
\xymatrix{
\ar@{}[r]|-{\In{\ ,\ }_v:}&H^1(F_v, \cG) \times H^1(F_v, \cG) \ar[r] \ar[d]^-{\wr}_-{\id \times H^1(F_v, \theta)} &H^2(F_v, \bG_m) \ar@{^(->}[r]^-{\inv_v} \ar@{=}[d] &\bQ/\bZ \ar@{=}[d]\\
&H^1(F_v, \cG) \times H^1(F_v, \cH) \ar[r]^-{\eqref{Shatz}} &H^2(F_v, \bG_m) \ar@{^(->}[r]^-{\inv_v} &\bQ/\bZ.
}
\ea\ee

\item \lab{maxl-isot}
One has $q_v(\Sel^j(\cG_{F_v})) = 0$ for $j \in \{ 1, 2\}$.

\item \lab{glob-quad}
For every $x \in H^1(U, \cG)$, its pullbacks $x_v \in H^1(F_v, \cG)$ satisfy $\sum_{v\not\in U} q_v(x_v) = 0$.
\eenum
\end{itemize}
Since $\Sel^j(\cG_{F_v})$ and $\Sel^j(\cH_{F_v})$ are orthogonal complements, these conditions ensure that for every $v \not \in U$ the subgroups $\Sel^j(\cG_{F_v}) \subset H^1(F_v, \cG)$ are maximal isotropic for the quadratic form $q_v$.
\epp

\brem \lab{quad}
In practice, the quadratic forms $q_v$ come into play only when $\#\cG$ is even. Indeed, for odd $\#\cG$, the equality $\In{x, x}_v = 2q_v(x)$ determines $q_v$, so if one insists that $\theta$ is such that the pairing
\[
b_{F_v}(-, \theta_{F_v}(\cdot))\colon \cG_{F_v} \times_{F_v} \cG_{F_v} \ra \bG_m \quad\quad\quad \text{for $v \not\in U$}
\]
is antisymmetric, then \ref{bil-pair} holds because the bilinear form $\In{\ ,\, }_v$ is symmetric (as may be seen using \v{C}ech cohomology), %(which in our setting agrees with the derived functor cohomology due to \cite{Sha64}*{Thm.~1}).
\ref{maxl-isot} follows from the isotropy of the subgroups $\Sel^j(\cG_{F_v})$, \ref{glob-quad} follows from the reciprocity for the Brauer group, whereas \ref{quad-form} results from the relation between $q_v$ and $\In{\ ,\, }_v$.
\erem

\bthm \lab{beyond-KMR}
In the setup of {\upshape\S\ref{break}}, 
\be\lab{beyond-KMR-eq}
\f{\# \Sel^1(\cG)}{\#\Sel^2(\cG)} \equiv \prod_{v\not\in U} \#\p{ \f{\Sel^1(\cG_{F_v})}{\Sel^1(\cG_{F_v}) \cap \Sel^2(\cG_{F_v})}  } \mod \bQ^{\times 2}.
\ee
\ethm

\bpf
The proof will use the following \Cref{XYZ}, which is a variant of \cite{KMR13}*{Lemma~2.3}. Before stating \Cref{XYZ}, we introduce the quadratic space
\[
(V, q) = \p{\textstyle\bigoplus_{v \not\in U} H^1(F_v, \cG),\, \sum_{v\not\in U} q_v}, \quad \text{whose associated bilinear form is} \quad \In{\ ,\, } \ce \textstyle\sum_{v\not\in U} \In{\ ,\, }_v.
\] 
Since $\In{\ ,\, }$ is continuous and nondegenerate, it exhibits $V$ as its own Pontryagin dual. The subgroups
\[
X = \textstyle\bigoplus_{v \not\in U} \Sel^1(\cG_{F_v}), \quad\quad\quad\quad Y = \textstyle\bigoplus_{v \not\in U} \Sel^2(\cG_{F_v}),\quad\quad\quad\quad Z = \im(\loc^1(\cG)),
\]
are maximal isotropic for $q$: indeed, $X$ and $Y$ due to \S\ref{break}~\ref{maxl-isot}, and $Z$ due to \Cref{im-orth} and \S\ref{break}~\ref{glob-quad}.

Since $\Sel^1(\cG)$ is finite (see \Cref{beyond-MR}), so is its quotient $X \cap Z$, and likewise for $Y \cap Z$. In $V$, both $X$ and $X\cap Y$ are open and compact, so $\f{X + Y}{Y} \cong \f{X}{X\cap Y}$ is finite. Thus, since 
\[
\tst \f{(X + Y)\cap Z}{Y \cap Z} \hra \f{X + Y}{Y},
\]
the intersection $(X + Y)\cap Z$ is finite, too.

\begin{lem} \lab{XYZ}
We have 
\[
\#\p{(X + Y)\cap Z} \equiv \#\p{(X\cap Z) + (Y \cap Z)} \bmod \bQ^{\times 2}.
\]
\end{lem}

\bpf
The proof of \cite{KMR13}*{Lemma 2.3} extends; we outline this extension.

The vanishing of the restrictions $\In{\ ,\, }|_X$ and $\In{\ ,\, }|_Y$ allows us to define a $\bQ/\bZ$-valued bilinear pairing $[\ ,\, ]$ on $(X + Y) \cap Z$ by
\[
[x + y, x\pr + y\pr] \ce \In{x, y\pr}\quad \text{for} \quad x + y,\, x\pr + y\pr \in (X + Y)\cap Z \quad  \text{with}\quad x, x\pr \in X \quad \text{and} \quad y, y\pr \in Y.
\]
The isotropy of $Z$, $X$, and $Y$ gives the vanishing
\[
\In{x, y} = q(x + y) - q(x) - q(y) = 0,
\] 
so $[\ ,\, ]$ is alternating. The resulting antisymmetry of $[\ ,\, ]$ ensures that the right and left kernels of $[\ ,\, ]$ agree; in particular, this common kernel $K$ contains $(X \cap Z) + (Y \cap Z)$. We claim that also
\be \lab{ker-rev}
K  \subset (X \cap Z) + (Y \cap Z), \qq \text{so that} \qq K = (X \cap Z) + (Y \cap Z).
\ee
To argue \eqref{ker-rev}, we fix $x$, $y$ as above with 
\[
x + y \in K, \qq \text{so that} \qq x \in (((X + Y)\cap Z) + (X \cap Y))^\perp,
\]
where the orthogonal complement is taken in $(V, \In{\ ,\,})$. Since the appearing subgroups are closed, \cite{HR79}*{24.10} gives
\[
(((X + Y)\cap Z) + (X \cap Y))^\perp = ((X \cap Y) + Z) \cap (X + Y) = ((X + Y)\cap Z) + (X \cap Y),
\]
so the freedom of adjusting $x$ and $y$ by opposite elements of $X \cap Y$ allows us to assume that 
\[
x\in (X +Y) \cap Z \subset Z.
\]
Then $y \in Z$ as well, which gives $x + y \in (X \cap Z) + (Y \cap Z)$, and \eqref{ker-rev} follows.

In conclusion, $[\ ,\, ]$ induces a nondegenerate alternating bilinear pairing on the quotient
\[
\tst \f{(X + Y)\cap Z}{(X \cap Z) + (Y \cap Z)}.
\] 
A well-known argument, as for example \cite{Dav10}*{proof of Lemma 4.2}, then implies that $\f{(X + Y)\cap Z}{(X \cap Z) + (Y \cap Z)}$ is a square of another finite abelian group, and hence is of square order, as desired. 
\epf

According to \Cref{beyond-MR}, the right hand side of \eqref{beyond-KMR-eq} equals $\#\p{\f{\Sel^{1 + 2}(\cG)}{\Sel^{1\cap 2}(\cG)}}$, which in turn equals $\#\p{\f{(X + Y) \cap Z}{X\cap Y\cap Z}}$. Moreover,
\[
\#\p{\f{(X + Y) \cap Z}{X\cap Y\cap Z}} \overset{\ref{XYZ}}{\equiv} \#((X\cap Z) + (Y \cap Z)) \cdot \#(X\cap Y \cap Z) \equiv \f{\#(X \cap Z)}{\#(Y \cap Z)} \mod \bQ^{\times 2},
\]
and it remains to observe that $\f{\#(X \cap Z)}{\#(Y \cap Z)} = \f{\#\Sel^1(\cG)}{\#\Sel^2(\cG)}$.
\epf

%\brem
%If $\Char F \nmid \#\cG$, then \cite{Mil06}*{II.2.13 (a)} implies the assumed finiteness of $\Sel^j(\cG)$.
%\erem

\section{Comparing Selmer sizes modulo squares in the main case of interest} \lab{ade}

In this section we specialize the results of \S\ref{ad-pre} to the setup of \Cref{main-eg}, which we assume and recall: $\cB$ and $\cB^\vee$ are global N\'{e}ron models of dual abelian varieties $B$ and $B^\vee$ that have good reduction at all the points of $U$, and 
\[
\cG = \cB[n]_U, \qq  \cH = \cB^\vee[n]_U
\]
for some $n \in \bZ_{> 0}$. Our main task is to justify that under suitable restrictions on $B$ various general assumptions made in \S\ref{ad-pre} are met if the local conditions are chosen as in \Cref{eg-loc-c} (we recall the choices in \Cref{1-2}). This justification leads to \Cref{dual}, which in \S\ref{final} will be a key ingredient in the proof of \Cref{main}.

\bprop \lab{1-2}
The following subgroups satisfy the hypotheses of {\upshape\S\ref{not-in-U}}, i.e.,~are open compact orthogonal complements.
\benum
\item \lab{1-2-a}
For $v\not\in U$, 
\[
\qq \Sel^1(\cG_{F_v}) = B(F_v)/nB(F_v) \quad \text{ and }\quad \Sel^1(\cH_{F_v}) = B^\vee(F_v)/nB^\vee(F_v).
\]
With these choices, 
\[
\qq \Sel^1(\cG) = \Sel_n(B) \qq \x{and} \qq \Sel^1(\cH) = \Sel_n (B^\vee);
\]
both Selmer groups are~finite.

\item \lab{1-2-b}
For $v\not\in U$,
\[\ba
\qq &\Sel^2(\cG_{F_v}) = H^1(\cO_v, \cB[n])  \quad\ \, \text{ and }\quad \Sel^2(\cH_{F_v}) = H^1(\cO_v, \cB^{\vee}[n])\quad \quad \quad\ \text{if $v \nmid \infty$}, \\
&\Sel^2(\cG_{F_v}) = B(F_v)/nB(F_v)  \quad \text{ and }\quad \Sel^2(\cH_{F_v}) = B^\vee(F_v)/nB^\vee(F_v)\quad\quad \text{if $v \mid \infty$}
\ea\]
under the assumption that $B$ {\upshape(}and hence also $B^\vee${\upshape)} has semiabelian reduction at all $v\in S\setminus U$ with $\Char \bF_v \mid n$. With these choices, $\Sel^2(\cG)$ is the subgroup $H^1(S, \cB[n])\pr \subset H^1(S, \cB[n])$ consisting of the elements whose restrictions to every $H^1(F_v, B[n])$ with $v\mid \infty$ lie in $B(F_v)/nB(F_v)$, and similarly for $\Sel^2(\cH)${\upshape;} both $\Sel^2(\cG)$ and $\Sel^2(\cH)$ are finite.
\eenum 
\eprop

\bpf 
The finiteness claims follow from the rest and from \Cref{beyond-MR}.
\benum
\item 
The orthogonal complement claim is a well-known important step of the proof of Tate local duality, compare \cite{Mil06}*{proof of III.7.8}. By \cite{Ces15b}*{4.2 and 4.3~(b)}, the connecting homomorphism
\[
\qq B(F_v) \ra H^1(F_v, B[n])
\]
is continuous and open (and likewise for $B^\vee$). Therefore, its image is open. Since $B(F_v)$ is compact, this image is also compact. The identification $\Sel^1(\cG) = \Sel_n B$ results from the definition of $\Sel_n B$ and from \cite{Ces16c}*{2.5 (d) and 4.2} (and similarly for $\Sel^1(\cH)$).

\item 
The $v\mid \infty$ case follows from \ref{1-2-a}, so we assume that $v\in S \setminus U$ (i.e.,~that $v\nmid \infty$). For such $v$, the given $\Sel^2$ are subgroups due to \cite{Ces16c}*{4.4} (which does not use the assumption on semiabelian reduction). Their openness and compactness follow from \cite{Ces15b}*{3.10}.

The orthogonal complement claim is essentially \cite{McC86}*{4.14}, but we must check that \eqref{Shatz} agrees with the pairing used in \emph{loc.~cit.} More precisely, the cohomology with supports sequence
\[
\qq \dotsb \ra H^m_{\bF_v}(\cO_v, \cB[n]) \ra H^m(\cO_v, \cB[n]) \ra H^m(F_v, B[n]) \ra H^{m + 1}_{\bF_v}(\cO_v, \cB[n]) \ra \dotsb
\]
of \cite{Mil06}*{III.0.3 (c)} is exact and the pairings in the diagram
\be\ba\lab{flat-ort}
\xymatrix{
H^2_{\bF_v}(\cO_v, \cB[n]) \ar@{}[r]|-{\bigtimes} & H^1(\cO_v, \cB^\vee[n]) \ar[d]\ar[rr]^-{\text{\cite{McC86}*{4.14}}} && H^3_{\bF_v}(\cO_v, \bG_m) \cong \bQ/\bZ \\
H^1(F_v, B[n]) \ar[u]\ar@{}[r]|-{\bigtimes} & H^1(F_v, B^\vee[n]) \ar[rr]^-{\eqref{Shatz}} && H^2(F_v, \bG_m) \cong \bQ/\bZ \ar@<15pt>[u]_-{\wr}  
}
\ea\ee
are perfect, so it suffices to prove that \eqref{flat-ort} and its analogue for $\cB^\vee[n]$ commute. Let 
\[
\qq i\colon \Spec \bF_v \hra \Spec \cO_v \qq \x{and} \qq j\colon \Spec F_v \hra \Spec \cO_v
\]
be the indicated immersions. Both pairings in \eqref{flat-ort} are Yoneda edge products---the top one due to its definition and the bottom one due to the observations made in the proof of \Cref{im-orth}---so the commutativity of the diagram \eqref{flat-ort} will follow from that of
\be\ba\lab{ext-fest}
\xymatrix{
\Ext^2(i_* \underline{\bZ}_{\bF_v} , \cB[n]) \ar@{}[r]|-{\bigtimes} & \Ext^1(\cB[n], \bG_m) \ar@{=}[d]\ar[rr]^-{\text{\cite{McC86}*{4.14}}} && \Ext^3(i_*\underline{\bZ}_{\bF_v}, \bG_m)  \\
\Ext^1(j_! \underline{\bZ}_{F_v} , \cB[n]) \ar[u]\ar@{}[r]|-{\bigtimes} & \Ext^1(\cB[n], \bG_m) \ar[d]\ar[rr] && \Ext^2(j_!\underline{\bZ}_{F_v}, \bG_m) \ar[u]_{\wr}  \\
\Ext^1(\underline{\bZ}_{F_v}, B[n]) \ar@{}[u]|-{\reflectbox{\rotatebox[origin=c]{90}{\text{\scalebox{1.5}{$\cong$}}}}}\ar@{}[r]|-{\bigtimes} & \Ext^1(B[n], \bG_m) \ar[rr]^-{\eqref{Shatz}} && \Ext^2(\underline{\bZ}_{F_v}, \bG_m), \ar@{}[u]|-{\reflectbox{\rotatebox[origin=c]{90}{\text{\scalebox{1.5}{ $\cong$}}
}}}
}
\ea\ee
where the identifications arise from the adjunction $j_! \dashv j^*$ as in \cite{Mil06}*{proof of III.0.3~(b)}. To see the commutativity of the bottom part of \eqref{ext-fest}, we replace $\cB[n]$ and $\bG_m$ by injective resolutions over $\Spec \cO_v$, interpret elements of $\Ext$ groups as homotopy classes of maps (compare with the proof of \Cref{im-orth} for this), and use the  adjunction $j_! \dashv j^*$ together with the fact that $j^*$ preserves injectives. 
To see the commutativity of the upper part, we observe that in the derived category the upper vertical arrows correspond to precomposition with the first morphism of the distinguished triangle 
\[
\qq i_* \underline{\bZ}_{\bF_v}[-1] \ra j_! \underline{\bZ}_{F_v} \ra \underline{\bZ}_{\cO_v} \ra i_* \underline{\bZ}_{\bF_v}.
\]

The identification 
\[
\qq \Sel^2(\cG) = H^1(S, \cB[n])\pr
\]
results from \cite{Ces16c}*{4.4}. 
\qedhere
\eenum
\epf

\brem \lab{rem-pr}
In many cases 
\[
H^1(S, \cB[n])\pr = H^1(S, \cB[n]);
\]
for instance, this happens if $n$ is odd or if $B(F_v)$ is connected for every real $v$ (which implies the same for $B^\vee$, cf.~\cite{GH81}*{\S1}) because then $H^1(F_v, B[n]) = 0$ for all $v\mid \infty$, as \emph{loc.~cit.}~proves. We resort to the somewhat artificial $H^1(S, \cB[n])\pr$ to make our duality results apply even when $H^1(F_v, B[n]) \neq 0$ for some $v\mid \infty$.
\erem

We turn to the assumptions of \S\ref{break} and note the following common source of suitable $\theta$.

\bprop \lab{anti-pf}
Let 
\[
\theta\colon \cB[n]_U \isomto \cB^\vee[n]_U
\]
be the isomorphism induced by a self-dual isogeny $\wt{\theta}$ of degree prime to $n$. Then $H^1(F_v, \theta)$ identifies the subgroup  $B(F_v)/nB(F_v)$ with $B^\vee(F_v)/nB^\vee(F_v)$ for all $v$ and identifies the subgroup $H^1(\cO_v, \cB[n])$ with $H^1(\cO_v, \cB^\vee[n])$ for all $v\nmid \infty$. \hfill\ensuremath{\qed}
\eprop

To address the remaining assumptions of \S\ref{break}, we first construct the quadratic forms $q_v$.

\bpp[Suitable quadratic forms $q_v$] \lab{qv-main}
Suppose that the assumptions of \Cref{1-2}~\ref{1-2-b} are met and that there is a self-dual isogeny 
\[
\wt{\theta}\pr \colon B \ra B^\vee
\]
of degree prime to $n$. Consider the self-dual isogeny 
\[
\wt{\theta} \ce \begin{cases} 2\wt{\theta}\pr, \text{ if $n$ is odd,} \\ \wt{\theta}\pr,\ \, \text{ if $n$ is even,} \end{cases}
\]
which also has degree prime to $n$. Due to \cite{PR12}*{Rem.~4.5}, the self-dual isogeny $\lambda \ce n\wt{\theta}$ comes from a symmetric line bundle $\sL$ on $B$, so the results of \cite{PR12}*{\S4} apply.  In particular, for $v \not \in U$, we can use the pullback $\sL_v$ of $\sL$ to $B_{F_v}$ to define the quadratic form 
\[
\xymatrix{
q_v\colon H^1(F_v, B[n]) \ar@{^(->}[r] &H^1(F_v, B[\gL]) \ar[r]^-{-\wt{q}_v} &H^2(F_v, \bG_m) \ar@{^(->}[r]^-{\inv_v} &\bQ/\bZ,
}
\]
where $\wt{q}_v\colon H^1(F_v, B[\gL]) \ra H^2(F_v, \bG_m)$ is the quadratic form provided by \cite{PR12}*{Cor.~4.7}.
\epp

\bprop \lab{qv-grand}
Assume the setup of {\upshape\S\ref{qv-main}}.
\benum
\item \lab{qv-grand-a}
The bilinear pairing associated to $q_v$ is the $\In{\ ,\, }_v$ of \eqref{pair-comp} with $\theta$ supplied by Proposition {\upshape\ref{anti-pf}}.

\item \lab{qv-grand-b}
We have $q_v(B(F_v)/nB(F_v)) = 0$ for every $v\not\in U$.

\item \lab{qv-grand-c}
We have $\sum_{v\not\in U} q_v(x_v) = 0$ for every $x \in H^1(U, \cB[n])$ with pullbacks $x_v \in H^1(F_v, B[n])$.

\item \lab{qv-grand-d}
We have $q_v(H^1(\cO_v, \cB[n])) = 0$ for every $v \in S \setminus U$ for which
\begin{enumerate}[label={(\roman*)}]
\item \lab{ass-i}
if $\Char \bF_v \mid n$, then $B$ has semiabelian reduction at $v$, and

\item \lab{ass-ii}
if $n$ is even, then the local Tamagawa factor $\#\Phi_v(\bF_v)$ is odd {\upshape(}here $\Phi_v = \cB_{\bF_v}/\cB_{\bF_v}^0${\upshape)}.
\eenum
\eenum
In particular, if {\upshape\ref{ass-i}} and {\upshape\ref{ass-ii}} hold for every $v\in S\setminus U$, then the $q_v$ meet the assumptions of {\upshape\S\ref{break}}.
\eprop

\bpf \hfill
\benum
\item 
By \cite{PR12}*{Cor.~4.7}, the bilinear pairing associated to $q_v$ is the restriction to $H^1(F_v, B[n])$ of the cup product pairing 
\[
\qq H^1(F_v, B[\gL]) \times H^1(F_v, B[\gL]) \ra \bQ/\bZ
\]
that uses Cartier duality 
\[
\qq b_{B[\gL]}\colon B[\gL] \times B[\gL] \ra \bG_m.
\]
Due to the naturality of the cup product, it remains to show that the diagram
\[
\xymatrix{
B[n] \times B[n] \ar[r]^-{\id \times \theta_F} \ar@{^(->}[d] & B[n] \times B^\vee[n] \ar@<-4pt>[d]^-{b_F} \\
B[\gL] \times B[\gL] \ar[r]^-{b_{B[\gL]}} & \bG_m
}
\]
commutes. For this, apply \cite{Oda69}*{Cor.~1.3 (ii)} with $\gA = \id_B$, $\gB = \wt{\theta}$, $\gL = [n]_B$, and $\gL\pr = \lambda$.

\item
This follows from \cite{PR12}*{Prop.~4.9} because the inclusion
\[
\qq H^1(F_v, B[n]) \hra H^1(F_v, B[\gL])
\] 
maps the subgroup $B(F_v)/nB(F_v)$ into $B^\vee(F_v)/\gL B(F_v)$ via the homomorphism induced by~$\theta$.

\item
By \cite{Ces16c}*{4.2 and 2.5 (d)}, for every 
\[
\qq x\in H^1(U, \cB[n]) \subset H^1(F, B[\gL])
\]
and every $v \in U$ one has 
\[
\qq x_v \in B(F_v)/nB(F_v) \subset B^\vee(F_v)/\gL B(F_v).
\]
Therefore, \cite{PR12}*{Thm.~4.14 (a)} gives the~claim.

\item
Set $n\pr \ce n$ if $n$ is odd, and $n\pr \ce \#\Phi_v(\bF_v)$ if $n$ is even. By \cite{Ces16c}*{2.5~(a)},
\[
n\pr x \in B(F_v)/nB(F_v)\quad \text{for every}\quad x \in H^1(\cO_v, \cB[n]).
\] 
Thus, \ref{qv-grand-b} gives $n^{\prime 2} q_v(x) = 0$. On the other hand, $\In{\ ,\,}_v$ vanishes on $H^1(\cO_v, \cB[n])$ due to \ref{qv-grand-a} and \S\ref{break} via \Cref{anti-pf}, so $2 q_v(x) = \In{x, x}_v = 0$. Since $n\pr$ is odd, we get $q_v(x) = 0$. \qedhere
\eenum
\epf

We are ready for the arithmetic duality result that will be used in the proof of \Cref{main}.

\bthm \lab{dual}
Fix an $n \in \bZ_{\ge 1}$, let $B$ be an abelian variety over a global field $F$, and let $\cB \ra S$ be its global N\'{e}ron model. For $v\nmid \infty$, let $\Phi_v$ denote the component group scheme of $\cB_{\bF_v}$. Suppose that 
\begin{enumerate}[label={{\upshape(\roman*)}}]
\item
$B$ has a self-dual isogeny $\wt{\theta}\pr$ of degree prime to $n$,

\item 
$B$ has semiabelian reduction at every nonarchimedean place $v$ of $F$ with $\Char \bF_v \mid n$, and

\item \lab{dual-iii}
if $n$ is even, then $\#\Phi_v(\bF_v)$ is odd for every nonarchimedean $v$.
\eenum
Let 
\[
H^1(S, \cB[n])\pr \subset H^1(S, \cB[n])
\]
be the subgroup of the elements whose restrictions to every $H^1(F_v, B[n])$ with $v\mid \infty$ lie in $B(F_v)/nB(F_v)$ {\upshape(}see Remark {\upshape\ref{rem-pr}} for some cases when $H^1(S, \cB[n])\pr = H^1(S, \cB[n])${\upshape)}.~Then
\[
\f{\#\Sel_n B}{\#H^1(S, \cB[n])\pr} \equiv \prod_{v\nmid \infty} \f{\#\Phi_v(\bF_v)}{\#(n\Phi_v)(\bF_v)} \bmod \bQ^{\times 2}.
\]
\ethm

\bpf
By \cite{Ces16c}*{2.5 (a)},
\[
\#\p{\f{B(F_v)/nB(F_v)}{H^1(\cO_v, \cB[n]) \cap (B(F_v)/nB(F_v))}} = \f{\#\Phi_v(\bF_v)}{\#(n\Phi_v)(\bF_v)} \qq \text{for $v \nmid \infty$}.
\]
Thus, the claim results from \Cref{beyond-KMR}, which applies due to \Cref{1-2,anti-pf,qv-grand}.
\epf

\brem \lab{l=2}
Since \Cref{beyond-KMR} is general, one would remove the assumption \ref{dual-iii} from \Cref{dual} by proving \Cref{qv-grand} \ref{qv-grand-d} without its assumption \ref{ass-ii}. This would also remove the additional assumption in the $\ell = 2$ case from \Cref{main}.
\erem

\section{The proof of \Cref{main}} \lab{final}

The goal of this section is to prove the following mild generalization of \Cref{main}.

\bthm \lab{main-pf}
Let $K$ be a global field of positive characteristic, let $\bF_q$ be its field of constants, let $\ell \nmid q$ be a prime, and let $A$ be an abelian variety over $K$. Suppose that $A_{K\ov{\bF}_q}$ has a polarization of degree prime to $\ell$ and, if $\ell = 2$, that the orders of the component groups of all the reductions of $A_{K\ov{\bF}_q}$ are odd. Then the $\ell$-parity conjecture holds for $A_{K\bF_{q^2}}$, i.e.,
\[
(-1)^{\rk_\ell(A_{K\bF_{q^2}})} = w(A_{K\bF_{q^2}}).
\]
\ethm

\bpf
We let $S$ be the connected smooth proper curve over $\bF_q$ having $K$ as its function field, we let $\cA \ra S$ and $\cA^\vee \ra S$ be the N\'{e}ron models of $A$ and $A^\vee$, and, for a variable closed point $v\in S$ (identified with the corresponding place $v$ of $K$), we let $\Phi_v$ denote the component group $\cA_{\bF_v}/\cA_{\bF_v}^0$ of the reduction $\cA_{\bF_v}$ of $A_{K_v}$. 

By \Cref{st-root}, for every even $n$ we have
\be \lab{root-even}
w(A_{K\bF_{q^n}}) = w(A_{K\bF_{q^2}}).
\ee
We claim that for every even $n$ we also have 
\be \lab{Sel-even}
\rk_\ell( A_{K\bF_{q^{n}}}) \equiv \rk_\ell (A_{K\bF_{q^2}}) \bmod 2.
\ee
Indeed, with the notation $\cX_\ell(-)$ of \S\ref{self-dual}, \Cref{DD} ensures that for every $1$-dimensional character $\chi$ of $\Gal(K\bF_{q^n}/K)$, the $\chi$-isotypical and the $\chi\i$-isotypical components of $\cX_\ell(A_{K\bF_{q^n}}) \tensor_{\bQ_\ell} \Qlbar$ have the same dimension. Therefore, the sum over all $\chi$ with $\chi^2 \neq 1$ of the $\chi$-isotypical components is even dimensional, whereas, by \Cref{sel-gal}, the sum over all $\chi$ with $\chi^2 = 1$ of such components is $\cX_\ell(A_{K\bF_{q^2}})\tensor_{\bQ_\ell} \Qlbar$. The claimed congruence \eqref{Sel-even} follows.
 
The combination of \eqref{root-even} and \eqref{Sel-even} allows us to replace $K$ by any $K\bF_{q^n}$ with $n$ even, so we loose no generality by assuming that $A$ has a polarization of degree prime to $\ell$, that the $\Gal(\ov{\bF}_q/\bF_q)$-action on $H^1_\et(S_{\ov{\bF}_q}, \cA[\ell])$ is trivial, and that $\Phi_v(\bF_v) = \Phi_v(\ov{\bF}_v)$ for every place $v$ of $K$.

By \Cref{rk-sel},
\be \lab{rk-exp} 
\rk_\ell (A_{K\bF_{q^2}}) \equiv \dim_{\bF_\ell}( \Sel_\ell (A_{K\bF_{q^2}})) - \dim_{\bF_\ell} (A[\ell](K\bF_{q^2}) )\bmod 2.
\ee
By \Cref{st-root},
\be \lab{root-exp}
w(A_{K\bF_{q^2}}) = (-1)^{\sum_{v\text{ inert in } K\bF_{q^2}} a(A_{K_v}) },
\ee
so our task becomes to compare the right sides of \eqref{rk-exp} and \eqref{root-exp}. We first use \Cref{cond-comp} to get
\be \lab{a-rewr}
\tst \sum_{v\text{ inert in } K\bF_{q^2}} a(A_{K_v}) = \sum_{v\text{ inert in } K\bF_{q^2}} a(A[\ell]_{K_{v}}) +  \sum_{v\text{ inert in } K\bF_{q^2}} \dim_{\bF_\ell} (\Phi_{v}[\ell](\bF_v)).
\ee
A place $v$ of $K$ splits in $K\bF_{q^2}$ if and only if the number of closed points of $S_{\ov{\bF}_q}$ above the closed point of $S$ determined by $v$ is even. Therefore, the Grothendieck--Ogg--Shafarevich formula \cite{Ray65}*{Thm.~1 (with (1 ter))} applied to $\cA[l]_{\ov{\bF}_q}$ (which is the N\'{e}ron model of its generic fiber, cf.~\cite{Ces16c}*{B.6}) gives the congruence
\be \lab{GOS-quot}
\tst \sum_{v\text{ inert in } K\bF_{q^2}} a(A[\ell]_{K_{v}}) \equiv \sum_{i = 0}^2  \dim_{\bF_\ell} (H^i_\et(S_{\ov{\bF}_q}, \cA[\ell])) \mod 2.
\ee
Let $j\colon U \hra S_{\ov{\bF}_q}$ be a nonempty open subscheme for which $\cA_U$ is an abelian scheme, so that $\cA[\ell]_U$ and $\cA^\vee[\ell]_U$ are Cartier dual by \cite{Oda69}*{Thm.~1.1}. The N\'{e}ron property ensures that 
\[
j_* (\cA[\ell]_U) \cong \cA[\ell]_{S_{\ov{\bF}_q}} \qq \x{and} \qq j_* (\cA^\vee[\ell]_U)  \cong \cA^\vee[\ell]_{S_{\ov{\bF}_q}}
\]
on the small \'{e}tale site of $S_{\ov{\bF}_q}$. Therefore, \cite{Mil80}*{V.2.2 (b)} supplies the duality isomorphism
\[
H^2(S_{\ov{\bF}_q}, \cA[\ell])^* \cong A^\vee[\ell](K\ov{\bF}_q).
\]
Since, in addition, a polarization of degree prime to $\ell$ gives an isomorphism 
\[
A^\vee[\ell](K\ov{\bF}_q) \simeq A[\ell](K\ov{\bF}_q),
\]
the congruence \eqref{GOS-quot} becomes
\be \lab{GOS-quot-2}
\tst \sum_{v\text{ inert in } K\bF_{q^2}} a(A[\ell]_{K_{v}}) \equiv \dim_{\bF_\ell} (H^1_\et(S_{\ov{\bF}_q}, \cA[\ell])) \mod 2.
\ee
We wish to express $\dim_{\bF_\ell} (H^1_\et(S_{\ov{\bF}_q}, \cA[\ell]))$ in terms of cohomology over $S_{\bF_{q^2}}$, so we use the Hochschild--Serre spectral sequence
\[
E_2^{ij} = H^i(\Gal(\ov{\bF}_q/\bF_{q^2}), H^j_\et(S_{\ov{\bF}_q}, \cA[\ell])) \Ra H^{i + j}_\et(S_{\bF_{q^2}}, \cA[\ell]),
\]
which degenerates on the $E_2$-page because $\bF_{q^2}$ has cohomological dimension $1$. The resulting short exact sequence
\[
0 \ra H^1(\Gal(\ov{\bF}_q/\bF_{q^2}), A[\ell](K\ov{\bF}_q)) \ra H^1_\et(S_{\bF_{q^2}}, \cA[\ell]) \ra H^1_\et(S_{\ov{\bF}_q}, \cA[\ell]) \ra 0
\]
together with the equalities
\[
\dim_{\bF_\ell} (H^1(\Gal(\ov{\bF}_q/\bF_{q^2}), A[\ell](K\ov{\bF}_q))) = \dim_{\bF_\ell} (H^0(\Gal(\ov{\bF}_q/\bF_{q^2}), A[\ell](K\ov{\bF}_q))) = \dim_{\bF_\ell} (A[\ell](K\bF_{q^2}))
\]
shows that
\[
\dim_{\bF_\ell} (H^1_\et(S_{\ov{\bF}_q}, \cA[\ell])) = \dim_{\bF_\ell} (H^1_\et(S_{\bF_{q^2}}, \cA[\ell])) - \dim_{\bF_\ell} (A[\ell](K\bF_{q^2})).
\]
Therefore, by combining this with \eqref{rk-exp}--\eqref{GOS-quot-2}, we find that our task is to prove that
\[
\tst \dim_{\bF_\ell} (\Sel_\ell (A_{K\bF_{q^2}})) \equiv \dim_{\bF_\ell} (H^1_\et(S_{\bF_{q^2}}, \cA[\ell])) + \sum_{v\text{ inert in } K\bF_{q^2}}  \dim_{\bF_\ell} (\Phi_{v}[\ell](\bF_v)) \mod 2.
\]
However, given our assumptions, this congruence is a special case of \Cref{dual}.
\epf

\end{document}